\documentclass[11pt]{amsart}
\usepackage{palatino}
\usepackage[all]{xy,xypic}

\usepackage{amsfonts,amsmath,oldgerm,amssymb,amscd}

\newcommand{\ra}{\rightarrow}		
\newcommand{\lra}{\longrightarrow}
\newcommand{\by}[1]{\stackrel{#1}{\ra}}

\newcommand{\remove}[1]{}

\newcommand{\ol}{\overline}		
\newcommand{\wt}{\widetilde}
\newcommand{\iso}{\by \sim}

\newtheorem{theorem}{Theorem}[section]
\newtheorem{proposition}[theorem]{Proposition}
\newtheorem{lemma}[theorem]{Lemma}
\newtheorem{definition}[theorem]{Definition}
\newtheorem{corollary}[theorem]{Corollary}
\newtheorem{conjecture}[theorem]{Conjecture}
\newtheorem{example}[theorem]{Example}

\newcommand{\ga}{\alpha}	\newcommand{\gb}{\beta}
		\newcommand{\gd}{\delta}

\newcommand{\BC}{\mbox{$\mathbb C$}}

	\newcommand{\mn}{\mbox{$\mathfrak n$}}

\newcommand{\ot}{\mbox{\,$\otimes$\,}}	
\newcommand{\op}{\mbox{$\oplus$}}
\newcommand{\Spec}{\text{Spec}}	
\newcommand{\hh}{\text{ht}}
\newcommand{\rank}{\text{rank}}

\newcommand{\dd}{\text{dim}}

\newcommand{\sur}{\twoheadrightarrow}
\newcommand{\bp}{\begin{proposition}}
\newcommand{\ep}{\end{proposition}}
\newcommand{\bl}{\begin{lemma}}
\newcommand{\el}{\end{lemma}}
\newcommand{\bt}{\begin{theorem}}
\newcommand{\et}{\end{theorem}}
\newcommand{\bc}{\begin{corollary}}
\newcommand{\ec}{\end{corollary}}
\newcommand{\bd}{\begin{definition}}
\newcommand{\ed}{\end{definition}}
\newcommand{\bco}{\begin{conjecture}}
\newcommand{\eco}{\end{conjecture}}

\renewcommand{\k}{\overline{\mathbb F}_p}

\def\rmk{\refstepcounter{theorem}\paragraph{{\bf Remark} \thetheorem}}
\def\proof{\paragraph{Proof}}

\def\notation{\paragraph{\bf Notation}}
\def\quest{\refstepcounter{theorem}\paragraph{{\bf Question} \thetheorem}}

\title[On two conjectures of Murthy]{On two conjectures of Murthy}
\author{ Mrinal Kanti Das}
\address{Stat-Math Unit, Indian Statistical Institute, 203 B. T. Road, Kolkata 700108 India}
\email{mrinal@isical.ac.in}
\date{\today}
\thanks{I sincerely thank Ravi Rao for some helpful comments and suggestions. }

\begin{document}

\begin{abstract}
This article concerns two conjectures of M. P. Murthy. For Murthy's conjecture on complete intersections, the major breakthrough  has still been the result proved by Mohan Kumar in 1978. In this article we improve \emph{Mohan Kumar's bound} when the base field is $\overline{\mathbb F}_p$, and illustrate some  applications of our result. Murthy's other conjecture is on a \emph{`splitting problem'}, which is roughly about
 finding the precise obstruction for a projective $R$-module $P$ of rank $\dd(R)-1$ to split off a free summand of rank one, where $R$ is a smooth affine algebra over an algebraically closed field $k$. Asok-Fasel achieved the initial breakthrough, by settling it for $3$-folds and $4$-folds when $char(k)\neq 2$. For $k=\k$ ($p\neq 2$) and $\dd(R)\geq 3$ we define an obstruction group and an obstruction class
for $P$ (whose determinant is trivial). As application we obtain: $P$ splits if and only if it maps onto a complete intersection ideal
of height $\dd(R)-1$.
\end{abstract}

\maketitle

\section{Introduction}
In this article we address two conjectures of M. P. Murthy. One of them is on the number of generators of an ideal in a polynomial ring over a field, and the other one is about splitting of a projective module defined over an affine algebra. We also study projective modules defined over certain threefolds. We discuss them one by one.

\medskip

\noindent
{\bf 1. Murthy's conjecture on complete intersections:} 
From a question of Murthy  in \cite{mu1}, the following conjecture has been formulated. Recall that $\mu(-)$ stands for the minimal number of generators.

\bco\label{mu1}(Murthy)
Let $k$ be a field and let $A=k[T_1,\cdots,T_d]$ be the polynomial ring in $d$ variables. Let $n\in {\mathbb N}$ and
$I\subset A$ be an ideal such that $\hh(I)=n=\mu(I/I^2)$. Then $\mu(I)=n$.
\eco

The conjecture is still open in general. 
The best known result on this conjecture  is due to Mohan Kumar \cite[Theorem 5]{mo}, stated below. It was proved in 1978.

\bt\cite{mo}\label{mohan1}
Let $k$ be a field and let $A=k[T_1,\cdots,T_d]$. Let $I\subset A$ be an ideal such that 
$\mu(I/I^2)=n\geq \dd(A/I)+2$. Then $\mu(I)=n$.
\et

Let $p$ be a prime  and $\k$ be the algebraic closure of the field consisting of $p$ elements. When $k=\k$, we improve Mohan Kumar's result in the following way (see  (\ref{monic}), (\ref{murthy}) in the text).

\bt\label{soln}
Let $I$ be an ideal in $A=\k[X_1,\cdots, X_d]$ with $\mu(I/I^2)=n\geq 2$. Assume that $n\geq  \dd (A/I)+1$. Then $\mu(I)=n$. 
\et

One can easily replace $\k$ by a PID $R$ which is of finite type over $\k$ in the above theorem. 
We now indicate one immediate application of (\ref{soln}).  

\bt\label{affine}
Let $\phi:\k[X_1,\cdots,X_5]\sur \k[Y_1,Y_2]$ be a surjective $\k$-algebra morphism and let 
$I$ be the kernel of $\phi$.  Then $I$ is generated by $3$ elements. 
\et

One may wonder whether $\k$ can be replaced by an arbitrary algebraically closed field in 
(\ref{affine}). It is in fact, a part of a famous open problem (see Section \ref{CI} and \cite{dug}).

\medskip

\noindent
{\bf 2. Murthy's splitting conjecture:} Based on a question asked by Murthy in \cite{mu2}, Asok-Fasel 
 formulated the following conjecture in \cite{af2} and termed it as \emph{Murthy's splitting conjecture}.

\bco\label{con2}(Murthy)
If $X=\Spec(R)$ is a smooth affine algebra of dimension $d$ over an algebraically closed field $k$ and $P$ is a projective $R$-module of rank $d-1$, then $P$ splits off a free summand of rank one (i.e., $P\iso Q\oplus R$ for some $R$-module $Q$) if and only if  the 
$(d-1)$\textsuperscript{th} Chern class $c_{d-1}(P)$  of $P$ vanishes in the Chow group $CH^{d-1}(X)$.
\eco

The case $d=2$ being trivial, we may focus on $d\geq 3$. The major breakthroughs on this conjecture have been achieved by Asok-Fasel for $d=3$ in \cite{af1}, and for $d=4$ in \cite{af2}, provided the characteristic of $k$ is unequal to $2$ in both the cases. Their methods involve ${\mathbb A}^{1}$-homotopy theory. The conjecture is open for $d\geq 5$.

According to the conjecture stated above, $CH^{d-1}(X)$ should serve as the obstruction group 
and $c_{d-1}(P)$ as the obstruction class for such splitting of $P$.
In this article we investigate such splitting  when $k=\k$, and $P$ has trivial determinant. 
Let  $X=\Spec(R)$ be a  smooth affine variety of dimension  $d\geq 3$ over $\k$ (with no restriction on $p$). Let $P$ be a projective $R$-module  of rank $d-1$, together with an isomorphism
$\chi:R\iso \wedge^{d-1}(P)$. In spirit of  Conjecture \ref{con2},
and following the footsteps of Bhatwadekar-Sridharan \cite{brs1, brs3,brs4},   we define a group $E^{d-1}(R)$ and associate to the pair $(P,\chi)$ an element $e_{d-1}(P,\chi)$ in $E^{d-1}(R)$ and prove the following result 
(see (\ref{split}) for a detailed statement):

\bt
The element $e_{d-1}(P,\chi)=0$ in $E^{d-1}(R)$ if and only if $P$ splits  a free summand of rank one.
\et

The group $E^{d-1}(R)$ is called the \emph{$(d-1)$\textsuperscript{th} Euler class group} of $R$ and $e_{d-1}(P,\chi)$
the \emph{$(d-1)$\textsuperscript{th} Euler class} of $(P,\chi)$. Our definition of $E^{d-1}(R)$ is built on  such a
 definition given in \cite{brs4}, with a vital modification. It will be interesting to explore the relation between the  groups $E^{d-1}(R)$ and $CH^{d-1}(X)$.

As an application of the Euler class theory developed in Section \ref{ECG}, combined with the cancellation theorem of Fasel-Rao-Swan \cite[7.5]{frs}, we are able to give the splitting problem a more tangible form by relating it to complete intersection ideals of height $d-1$ in $R$. More precisely, we prove the following theorem (see  (\ref{murthy2})). Note that this answers a weaker form of (\ref{con2}).

\bt\label{main}
Let $p\neq 2$, and $R$ be a smooth  affine algebra of dimension $d\geq 3$ over $\k$, 
and $P$ be a projective $R$-module of rank $d-1$ with trivial determinant.  Then $P\iso Q\op R$ for some $R$-module $Q$ if and only if there
is a surjection $\alpha:P\sur J$ such that $J\subset R$ is an ideal of height $d-1$ and
 $\mu(J)=d-1$.
\et


We use the smoothness of $R$ only to apply the $3$-fold case of \cite[7.5]{frs}. 
However, Ravi Rao has recently informed us (private communication) that the smoothness assumption 
in \emph{op. cit.} has since been removed. Here in this article, we are able to show that we do not need any smoothness assumption in (\ref{main}) if $d\geq 5$ and $gcd(p,(d-1)!)=1$
(see (\ref{murthy2a}) in the text).

\medskip

\noindent
{\bf 3. A question on threefolds:} Let $k$ be an algebraically closed field and $R$ be an affine $k$-algebra of dimension $3$. Let $P$ be a projective $R$-module of rank $2$
with trivial determinant. By a theorem of Eisenbud-Evans \cite{ee}, there is a surjection $\alpha:P\sur J$, where $J\subset R$ is an ideal of height $2$. We shall call such $J$ to be a \emph{generic section} of $P$. Obviously, any $R$-module in the isomorphism class of $P$ has the same generic sections as $P$.
In this context, we may pose the following question.

\smallskip

\quest\label{quest}
\emph{Does a generic section  of $P$ uniquely determine the isomorphism class of $P$? In other words, if
$Q$ is another projective $R$-module of rank two with trivial determinant such that $P$ and $Q$ have
a common generic section  $J$, then is $P$ isomorphic to $Q$?}

\medskip

We answer this question affirmatively when $k=\k$, $p\neq 2$,  and $R$ is smooth. Again, smoothness might be a superfluous assumption.
We would also like to point out to the reader that apart from \cite{frs}, we have not used any other cancellation results (e.g., those due to Asok-Fasel  from \cite{af1}). On the other hand, as an application of our results we easily derive that the projective $R$-modules of rank $2$ with trivial determinant are cancellative, where $R$ is a smooth affine $3$-fold over $\k$ and $p\neq 2$.

As it turns out, as opposed to \cite{af1,af2}, we need much lighter machinery to prove the results in this paper. It is unlikely that the methods of this paper can be extended to the case when $\k$ is replaced by an arbitrary algebraically closed field. The great advantages of working with $\k$ are the following: \emph{(1) For a
one dimensional affine $\k$-algebra $A$, the group $SK_1(A)$ is trivial by a result due to Krishna-Srinivas \cite{ks}.   (2) For a two-dimensional affine $\k$-algebra $A$, any projective $A$-module with trivial determinant is free (\ref{free}).}
We exploit these two results heavily in this article. Note that by \cite[2.5]{mup}, (1) is no longer valid if $\k$ is replaced by $\BC$, whereas (2) is well-known to be false for 
affine surfaces over $\BC$.

\section{Some assorted results}
The purpose of this section is  to collect various results from the literature. 
Quite often we would tailor them or improve them a little bit to suit our requirements in subsequent sections. We start with a result from \cite{ks} which is crucial to this paper.

\bt\label{sk1}
Let $R$ be an affine algebra of dimension one  over $\k$. Then, $SK_1(R)$ is trivial.
\et

\proof
See the last two paragraphs of the proof of \cite[Theorem 6.4.1, page 274]{ks}.
\qed

\smallskip

The next result is due to Swan.

\bp\cite[9.10]{sw}\label{swan}
Let $A$ be a ring and $I$ be an ideal. Let $\gamma\in Sp_{2t}(A/I)$, $t\geq 1$. If the class of 
$\gamma$ is trivial in $K_1Sp(A/I)$ and if  $2t\geq sr(A)-1$ (where $sr(-)$ means stable range),
then $\gamma$ has a lift $\ga\in Sp_{2t}(A)$. 
\ep

The following corollary (of (\ref{sk1}) and (\ref{swan}))  must be well-known but  we did not
find any suitable reference. 

\bc\label{vas}
Let $A$ be an affine algebra    over $\k$, and let $I\subset A$ be an ideal such that 
$\dd(A/I)\leq 1$. Then, we have the following assertions.
\begin{enumerate}
\item
The canonical map $SL_{n}(A)\lra SL_n(A/I)$ is surjective for $n\geq 3$.
\item
If $\dd(A)=3$, then the canonical map $SL_{2}(A)\lra SL_2(A/I)$ is surjective.
\end{enumerate}
\ec

\proof
If $\dd(A/I)=0$, then  $SL_n(A/I)=E_n(A/I)$, and we are done because $E_n(A)\lra E_n(A/I)$ is surjective for  $n\geq 2$.
Therefore, we assume that $\dd(A/I)=1$. 

Let $n\geq 3$. Then, applying  (\ref{sk1}) on $A/I$, together with  Vaserstein's stability results 
\cite{v}, we have $SL_{n}(A/I)=E_n(A/I)$ for $n\geq 3$ and we have thus proved (1). 

To prove (2), we need some additional arguments. Let $\dd(A)=3$, $\dd(A/I)=1$, and let
$\gamma\in SL_{2}(A/I)$ be arbitrary.
By the remark following \cite[16.2]{sv}, we have $K_1Sp(A/I)=SK_1(A/I)$ and therefore, by 
(\ref{sk1}),  $K_1Sp(A/I)$ is trivial.  Note that $\gamma\in SL_2(A/I)=Sp_2(A/I)$. By 
\cite[17.3]{sv}, we have $sr(A)\leq \text{ max }\{2,\dd(A)\}=3$. We can now
apply (\ref{swan}) with $t=1$ to obtain  $\alpha\in Sp_{2}(A)(=SL_{2}(A))$
which is a lift of $\gamma$. This 
completes the proof.
\qed

\smallskip

The next theorem is essentially an accumulation of results of various authors.

\bt\label{gen}
Let $R$ be an affine algebra of dimension $d\geq 2$ over $\k$ and $I\subset R$ be an ideal such that 
$\mu(I/I^2)=d$. Suppose it is given that $I=(a_1,\cdots,a_d)+I^2$. Then,
there exist $b_1,\cdots,b_d\in I$ such that $I=(b_1,\cdots,b_d)$ with $b_i-a_i\in I^2$ for 
$i=1,\cdots, d$.
\et

\proof
This result was proved in \cite{mamu} when $R$ is reduced and in addition, $R$ is smooth 
if $d=2$.
Let us first  assume that $R$ is reduced. In this case, we need only remove the smoothness assumption
when $d=2$. But it has been proved in \cite{ks} that $F^2K_0(R)$ is trivial even when $R$ is singular and therefore we can follow the proof of \cite{mamu}.

It is not difficult to prove that we can take  $R$ to be reduced to start with.
To see this, let $\mn$ be the nilradical of $R$
and let bar denote reduction modulo $\mn$. Note that,
$$\overline{I}=(\ol{a}_1,\cdots,\ol{a}_d)+\ol{I}^2 \text{ in } \ol{R}.$$ 
By the above paragraph, there exist $\ol{c}_1,\cdots,\ol{c}_d$ such that 
$\ol{I}=(\ol{c}_1,\cdots,\ol{c}_d)$ with $\ol{c}_i-\ol{a}_i\in \ol{I}^2$. One can now follow
the proof of \cite[4.13]{ke} (the `injectivity' part) 
to see that there exist $b_1,\cdots,b_d\in I$ such that $I=(b_1,\cdots,b_d)$ with $b_i-a_i\in I^2$
 for 
$i=1,\cdots, d$.
\qed

\medskip

Some immediate corollaries are in order. For the definition of the `top' Euler class group, see \cite{brs3}. Here we abbreviate the $d$\textsuperscript{\,th} Euler class group $E^d(R,R)$ as $E^d(R)$.  From
(\ref{gen}) and the definition of $E^d(R)$, we have:

\bc\label{ecg}
Let $R$ be an affine algebra of dimension $d\geq 2$ over $\k$. The $d$\textsuperscript{\,th} Euler class group 
$E^d(R)$ is trivial.
\ec

Let us recall the definition of \emph{reduction} of an ideal.

\smallskip

\bd
Let $A$ be a ring . Let $J\subseteq I$ be ideals. $J$ is said to be a {\it reduction} of $I$ if there exists a non-negative integer $t$ such that $I^{t+1}=JI^{t}$.
\ed

\bc
Let $R$ be an affine algebra of dimension $d\geq 2$ over $\k$ and $I\subset R$ be an ideal of height 
$d$. Then $I$ is a set-theoretic complete intersection. 
\ec

\proof
It follows from a theorem of Katz \cite{ka} (see the proof of \cite[8.73
(2)]{hs})  that $I$ has a reduction $J$ such that $J/J^2$ is generated by $d$ elements. Since  $J$ is a reduction of $I$, we have 
$\sqrt{J}=\sqrt{I}$ and $\text{ht}(J)=\text{ht}(I)$. By (\ref{gen}), $J$ is generated by $d$
elements. Therefore, $I$ is set-theoretically  generated by $d$
elements. 
\qed

\bc\label{split1}
Let  $R$ be an affine algebra of dimension $d\geq 2$ over $\k$ and $P$ be a projective $R$-module of rank $d$ with trivial determinant. Then $P\iso Q\oplus R$ for some $R$-module $Q$.
\ec

\proof
For $d\geq 3$, this is proved in \cite{mmr}. Let $d=2$. In this case, $R$ is assumed to be smooth in \cite{mmr}. We can  remove this assumption, as follows.

It is easy to see following the proof of \cite[3.1]{brs3} that for \emph{any} two-dimensional Noetherian ring 
$A$ and a projective $A$-module $M$ with trivialization $\chi:A\iso \wedge^2(M)$, the Euler class $e(M,\chi)$ is a well-defined element of $E^2(A)$. Recall from \cite[4.4]{brs3} that $e(M,\chi)=0$ if and only if $M$ splits a free summand.

By (\ref{ecg}), $E^2(R)$ is trivial and therefore, the proof is complete by 
\cite[4.4]{brs3}.
\qed

\smallskip

We shall see repeated use of the following corollary in this article.

\bc\label{free}
Let  $R$ be a two-dimensional affine algebra over $\k$, and $P$ be a projective $R$-module of rank $2$ with trivial determinant. Then $P$ is free. In fact, any projective $R$-module  with trivial determinant is free.
\ec

\proof
Let $\rank(P)=n\geq 2$.
By a classical result of Serre \cite{se}, $P\iso P'\op R^{n-2}$ for some projective $R$-module $P'$ of rank $2$. Note that the determinant of $P'$ is trivial. Now, $P'$ is free by  (\ref{split1}).
\qed

\remove{
\bc\label{ravi}
Let $p\neq 2$.  Let  $R$ be a two-dimensional affine algebra over $\k$. Then, any unimodular row 
$(f_1(T),f_2(T),f_3(T))$  over $R[T]$ can be completed to an invertible matrix.
\ec

\proof
By a result of Rao \cite{ra2}, the unimodular row $(f_1(T),f_2(T),f_3(T))$ is extended from $R$.
A unimodular row of length three over $R$ defines a stably free $R$-module $P$ of rank $2$.
Now $P$ is free by the above corollary.
\qed}

\smallskip

The following is a standard result. For a proof, the reader may see \cite[2.17]{ke}.

\begin{lemma}\label{26}
Let $A$ be a Noetherian ring and $P$ a finitely generated projective 
$A$-module. Let $P[T]$ denote projective $A[T]$-module $P\otimes A[T]$.
Let $\alpha (T) : P[T]\twoheadrightarrow A[T]$ and $\beta (T) : P[T]\twoheadrightarrow A[T]$
be two surjections such that $\alpha (0) = \beta (0)$. Suppose further that 
the projective $A[T]$-modules {\rm ker}\,$\alpha (T)$ and 
{\rm ker}\,$\beta (T)$ are extended from $A$. Then there exists an automorphism $\sigma (T)$ of $P[T]$ with
$\sigma (0)=$ id such that $\beta (T)\sigma (T) = \alpha (T)$.
\end{lemma}

The next lemma follows from the well known Quillen Splitting Lemma 
\cite[Lemma 1]{q}, and its proof is essentially contained in 
\cite[Theorem 1]{q}.

\bl\label{qsplit}
Let $A$ be a Noetherian ring and $P$ be a finitely generated projective 
$A$-module. Let $s,t\in A$ be such that $As+At=A$. Let $\sigma(T)$ be an
$A_{st}[T]$-automorphism of $P_{st}[T]$ such that $\sigma (0)=$ id.
Then, $\sigma(T) = {\alpha(T)}_{s}{\beta(T)}_{t}$, where
$\alpha(T)$ is an $A_{t}[T]$-automorphism of $P_{t}[T]$ such that
$\alpha(T) =$ id modulo the ideal $(sT)$ and 
$\beta (T) $ is an $A_{s}[T]$-automorphism of $P_{s}[T]$ such that
$\beta (T)=$ id modulo the ideal $(tT)$.
\el

We shall need the following \emph{``moving lemma"}. The version given below can easily be proved following \cite[2.4]{brs4}, which in turn is essentially  based on \cite[2.14]{brs3}.

\bl\label{move}
Let $A$ be a Noetherian ring of dimension $d$ and let $J\subset A$ be an ideal of height $n$
such that $2n\geq d+1$. Let $P$ be a projective $A$-module of rank $n$ and $\ol{\ga}:P/JP\sur J/J^2$ be a surjection. Then, there exists an ideal $J'$ of $A$ and a surjection $\gb:P\sur J\cap J'$ such that:
\begin{enumerate}
\item
$J+J'=A$,
\item
$\gb \ot A/J=\ol{\ga}$,
\item
$\hh(J')\geq n$.
\end{enumerate}
Given any ideal $K\subset A$ of height $n$, the map $\gb$ can be  chosen so that $J'+K=A$.  
\el

We need the following result from \cite{brs4}.

\bl\cite[5.1]{brs4}\label{hom}
Let $A$ be a Noetherian ring and $P$ a projective $A$-module of rank $n$. Let $\alpha:P\sur J$ and $\beta:P\sur J'$ be two surjections, where $J,J'$ are ideals of $A$ of height at least $n$. Then, there exists an ideal $I\subset A[T]$ of height $\geq n$ and a surjection $\phi(T):P[T]\sur I$ such that $I(0)=J$,
$\phi(0)=\ga$ and $I(1)=J'$,  $\phi(1)=\gb$. 
\el

Next we state a remarkable result of Rao \cite[3.1]{ra3}. Rao proved it when $A$ is local. The following version can be deduced from \cite{ra3} by applying Quillen's local-global principle \cite{q}. 

\bt\label{ravi}
Let $A$ be a Noetherian ring of dimension $3$  such that $2A=A$. Take any  unimodular row
$(f_1(T),f_2(T),f_3(T))$  over $A[T]$. Then there is $\theta(T)\in GL_{3}(A[T])$ such that 
$(f_1(T),f_2(T),f_3(T))\theta(T)=(f_1(0),f_2(0),f_3(0))$. Replacing $\theta(T)$ by 
$\theta(T)\theta(0)^{-1}$ we can actually conclude that $\theta(T)\in SL_3(A[T])$ and 
$\theta(0)=id$.
\et

The above theorem will enable us to cover the case $d=4$ in Sections \ref{AddSub} and 
\ref{ECG} as it facilitates certain patching arguments. We illustrate one such instance in
the following theorem, which  is a variant of \cite[5.2]{brs4}.

\bt\label{chi}
Let $A$ be a Noetherian ring of dimension $d\geq 3$ with the following additional assumptions: 
(i) $2A=A$ if $d=4$, (ii) $A$ is regular containing a field if $d\geq 5$.  Let $P$ be a projective $A$-module of rank $d-1$ and let  $\ga(T):P[T]\sur I$ be a surjection where $I$ is an ideal of height $d-1$. Assume that 
$J=I(0)$ is a proper ideal, and further that $P/NP$ is free,  where $N=(I\cap A)^2$. Let $p_1,\cdots,p_{d-1}\in P$ be such that their images in $P/NP$ form a basis. Let $a_1,\cdots,a_{d-1}\in J$ be such that 
$\ga(0)(p_i)=a_i$. Then, there exists an ideal $K\subset A$ of height $\geq d-1$ such that $K+N=A$ and:

\begin{enumerate}
\item
$I\cap K[T]=(F_1(T),\cdots,F_{d-1}(T))$,
\item
$F_i(0)-F_i(1)\in K^2$, $i=1,\cdots,d-1$,
\item
$\ga(T)(p_i)-F_i(T)\in I^2$,
\item
$F_i(0)-a_i\in J^2$.
\end{enumerate}
\et

\proof
This is essentially contained in \cite{brs4}, with $A$ regular (containing a field).  
Therefore, for $d\geq 5$, we are done. For the remaining cases, retaining the notations from \cite{brs4} we give a sketch of the necessary modification in \cite[page 151, second paragraph]{brs4}.

We have (from the proof of \cite[5.2]{brs4}), 
$I'=I\cap K[T]$. An element $a\in N$ is obtained so that $I'_{1+a}=I_{1+a}=(G_1(T),\cdots,
G_{d-1}(T))$. On the other hand, $I'_a=KA_a[T]=(a_1+cb_1,\cdots,a_{d-1}+cb_{d-1})$ such that 
$G_i(0)=a_i+cb_i$ for $i=1,\cdots,d-1$.

We now split the cases.

\smallskip

\noindent
{\it Case 1.} $d=3$. Let $b=1+a$. The rows $ (G_1(T),G_{2}(T))$ and $(a_1+cb_1,a_{2}+cb_{2})$ are unimodular over the ring 
$A_{ab}[T]$, and they agree when $T$ is set to zero.
As any unimodular row of length two over any ring can be completed to a $2\times 2$ with determinant one, we can find $\theta(T)\in SL_2(A_{ab}[T])$ such that $(G_1(T),G_2(T))\theta(T)=(1,0)$. Then $(G_1(0),G_2(0))\theta(0)=(1,0)$, implying
that $(a_1+cb_1,a_2+cb_2)\theta(0)=(1,0)$. Taking $\sigma(T):=\theta(T)\theta(0)^{-1}$, we
observe that $\sigma(0)=id$, and $(G_1(T),G_2(T))\sigma(T)=(a_1+cb_1,a_2+cb_2)=(G_1(0),G_2(0))$. The rest of the arguments are the same as \cite[5.2]{brs4}.

\smallskip

\noindent
{\it Case 2.} $d=4$ and $2A=A$.  Consider the unimodular rows $(G_1(T),G_2(T),G_{3}(T))$
and $(a_1+cb_1,a_2+cb_2,a_{3}+cb_{3})$ over the ring $A_{a(1+aA)}[T]$. Note that $\dd(A_{a(1+aA)})\leq 3$. By  (\ref{ravi}) there exists $\sigma'(T)\in SL_3(A_{a(1+aA)}[T])$ such that $\sigma'(0)=id$ and 
$$(G_1(T),G_2(T),G_{3}(T))\sigma'(T)=(a_1+cb_1,a_2+cb_2,a_{3}+cb_{3}).$$
 We can find  some $b$ of the form $1+\lambda a$ such that $b$ is a multiple of $1+a$, and 
some $\sigma(T)\in SL_3(A_{ab}[T])$ such that $\sigma(0)=id$ and $(G_1(T),G_2(T),G_{3}(T))\sigma(T)=
(a_1+cb_1,a_2+cb_2,a_{3}+cb_{3})$ over the ring $A_{ab}[T]$. The rest is same
as \cite[5.2]{brs4}.
\qed

\smallskip

\rmk
The assumption $2A=A$ is not required for $d=4$ if $A$ is regular containing a field.

\section{Results on Conjecture \ref{mu1} and their applications}\label{CI}
We start with a lemma of Mohan Kumar \cite{mo1}, recast slightly to suit our needs.

\bl\label{nak}
Let $A$ be commutative Noetherian ring and $J\subset A$ be an ideal. Assume that $J=K+L$, where
$K, L$ are ideals and $L\subset J^2$. Then, there exist $e\in L$ such that:
\begin{enumerate}
\item
$J=(K,e)$ and $e(1-e)\in K$,
\item
If $J'=(K,1-e)$, then $J\cap J'=K$ (note also  that $J'+L=A$),
\item
For any $a\in A$, we have $(J,a)=(K, e+(1-e)a)$.
\end{enumerate} 
\el

We first prove a result on ideals containing monic polynomials. Compare it with the results of Mohan Kumar and Mandal \cite{mo,m1}.
I learned slightly simplified proofs of their results from Raja Sridharan during my graduate days, and I mimic those arguments here. 

\bt\label{monic}
Let $R$ be an affine algebra over $\k$ and $n\geq 2$ be an integer. Let $I\subset R[T]$ be an ideal containing a monic polynomial  such that $\mu(I/I^2)=n\geq \dd (R[T]/I)+1$. Then $I$ is generated by $n$ elements. 
Moreover, any set of $n$ generators of $I/I^2$ can be lifted to a set of $n$ generators of $I$ in the following cases:
\begin{enumerate}
\item
 $n=2=\mu(I/I^2)=\hh(I)=\dd(R)$ and $p\neq 2$, 
\item
$n\geq 3$.
\end{enumerate} 
\et

\proof
Let $f\in I$ be a monic polynomial and assume that $I=(f_1,\cdots,f_n)+I^2$. We treat two cases separately.

\noindent
\emph{Case 1.} Let $n=2$. We have $I=(f_1,f_2)+I^2$. By Lemma \ref{nak} there exists $e\in I^2$ such that $e(1-e)\in (f_1,f_2)$. Therefore, $I_{1-e}=(f_1,f_2)_{1-e}$. On the other hand, $I_e=R[T]_e=(1,0)$. 
Note that the unimodular row $(f_1,f_2)_{e(1-e)}$ over $R[T]_{e(1-e)}$  can be completed to a $2\times 2$ matrix with determinant $1$.  Using a standard patching argument we obtain a projective $R[T]$-module $P$ of rank $2$ with \emph{trivial determinant} such that $P$ maps onto $I$. Since $I$ contains the monic polynomial $f$, the module $P_f$ has a unimodular element and hence becomes free. The Affine Horrocks Theorem implies that $P$ is free. Consequently, $I$ is generated by $2$ elements, say, $I=(g_1,g_2)$.

We now show that if $p\neq 2$, then  we can actually find $f'_1,f'_2\in I$ such that $I=(f'_1,f'_2)$ with 
$f'_i-f_i\in I^2$, $i=1,2$. We have, $I=(f_1,f_2)+I^2$ and
$I=(g_1,g_2)$ thus far. By \cite[2.2]{b}, there exists a matrix $\ol{\sigma}\in GL_{2}(R[T]/I)$
such that $(\ol{f}_1,\ol{f}_2)=(\ol{g}_1,\ol{g}_2)\ol{\sigma}$.  Let $\det(\ol{\sigma})=\ol{u}$. Let $v\in R[T]$ be such that $uv-1\in I$. Then $(v,g_2,-g_1)$ is a unimodular row over $R[T]$ and since $p\neq 2$, this row is completable by \cite[2.5]{ra2}. By \cite[2.3]{rs}, there are generators $h_1,h_2$ of $I$
and $\ol{\delta}\in GL_{2}(R[T]/I)$ such that $\det(\ol{\delta})=\ol{u}$ and $(\ol{h}_1,\ol{h}_2)=(\ol{g}_1,\ol{g}_2)\ol{\delta}$. Consequently, $(\ol{f}_1,\ol{f}_2)=(\ol{h}_1,\ol{h}_2)\ol{\tau}$ for some $\ol{\tau}\in SL_2(R[T]/I)$. We can now apply (\ref{vas}) and lift $\ol{\tau}$ to some $\tau\in SL_2(R[T])$. Write $(h_1,h_2)\tau=(f'_1,f'_2)$. Then $f'_1,f'_2$ is the desired set of generators for $I$.

\medskip

\noindent
\emph{Case 2.} Let $n\geq 3$. Adding $f^t$ to $f_1$ for suitable $t$
 we can assume that $f_1$ is monic. Let $J=I\cap R$ and $B=R[T]/(J^2[T],f_1)$. Then,
$$\overline{I}=(\overline{f}_2,\cdots,\overline{f}_n)+\overline{I}^2,$$
where bar means modulo $(J^2[T],f_1)$. Note that $B$ is an affine algebra over $\k$ and as 
$B\hookrightarrow R[T]/I$ is integral, we have $\dd (B)=\dd (R[T]/I)\leq n-1$. By (\ref{gen}), there
exist $h_2,\cdots,h_n\in I$ such that $\overline{I}=(\overline{h}_2,\cdots,\overline{h}_n)$. Then,
$$I=(f_1,h_2,\cdots,h_n)+J^2[T].$$
Applying (\ref{nak}) we find an ideal $I'\subset R[T]$ such that $I\cap I'= (f_1,h_2,\cdots,h_n)$ with 
$I'+J^2[T]=R[T]$. As $I'$ contains a monic (namely, $f_1$), by \cite[Lemma 1.1, Chapter III]{la}, $I'\cap R+J^2=R$. Let $s\in J$ 
be such that $1+s\in I'$. Then, 
$$I_{1+s}=(f_1,h_2,\cdots,h_n)_{1+s}.$$ 
Let $\phi:R_{1+s}[T]^n\sur I_{1+s}$ be the corresponding surjection. On the other hand, we have 
$I_{s}=R_{s}[T]$ and we have the surjection $\psi:R_s[T]^n\sur I_s=R_s[T]$ which sends the first basis 
vector to $1$ and the rest to zero. Now, over $R_{s(1+s)}[T]$ we have the unimodular row $(f_1,h_2,\cdots,h_n)$ whose one entry is monic. Since this row is elementarily completable by \cite{ra}, a standard patching argument completes the proof.
\qed

\medskip

Some immediate applications follow. We first record the following improvement on Mohan Kumar's result on Murthy's conjecture.

\bt\label{murthy}
Let $I$ be an ideal in $A=\k[X_1,\cdots, X_d]$ with $\mu(I/I^2)=n\geq \dd (A/I)+1$. Then $\mu(I)=n$.
\et

\proof
After a change of variables, $I$ contains a monic polynomial in one of the variables, say, 
$X_d$. Write $R=\k[X_1,\cdots, X_{d-1}]$ and apply the above theorem.
\qed

\medskip

\rmk
In (\ref{murthy}) we can replace $\k$ by a PID which is of finite type  over $\k$. 

\medskip

\rmk\label{murthy1}
In the above theorem, if  $n\geq 3$, then any given  set of $n$ generators of $I/I^2$ can be lifted to a set of $n$ generators of $I$. If $p\neq 2$, and $n=2=\mu(I/I^2)=\hh(I)=\dd(R)$,  then also we can 
lift generators of $I/I^2$ to generators of $I$ as in (\ref{monic}). In contrast,  
Bhatwadekar-Sridharan showed in
\cite[3.15]{brs1} that for the complete intersection ideal $I=(X_1^3+X_2^3-1,X_3-1)\subset \BC[X_1,X_2,X_3]$, there is a set of $2$ generators of $I/I^2$  which cannot be lifted to a set of two generators of $I$.  The crux of the matter in their example is that 
$SK_1(\BC[X_1,X_2]/(X_1^3+X_2^3-1))$ is non-trivial by \cite[2.5]{mup}.
\medskip

\begin{example}\label{sri}
In this example we show that  assertion (1) in (\ref{monic}) for $n=2$ (or a similar assertion in 
(\ref{murthy}) as described in (\ref{murthy1})) is not valid if $\dd(R)\geq 3$. In this sense, these results cannot be improved. The idea of this example is the same as \cite[3.15]{brs1}. Let $k$ be any algebraically closed field of characteristic unequal to $2$. Let
$$B=\frac{k[X,Y,Z]}{(Z^2-XY)}.$$
Srinivas \cite{sr} proved that $SK_1(B)\neq 0$   (In fact, he gave
example of  
an explicit matrix in $SL_2(B)$ which is not even in the image of $SL_2(k[X,Y,Z])\lra SL_2(B)$. His example is related to Cohn's example \cite{c}). Write $f=Z^2-XY$, $R=k[X,Y,Z]$, and consider the ideal  
$I=(f,T-1)\subset R[T]$, where $T$ is an indeterminate. Let bar denote reduction modulo $I$.
We have, $R[T]/I=B$.
By the aforementioned result of Srinivas \cite{sr}, there is a $\ol{\sigma}\in SL_2(R[T]/I)$ such that it is not stably elementary.
Let $(\ol{f},\ol{T-1})\ol{\sigma}=(\ol{G},\ol{H})$. Then, (i) $I$ is an ideal of height $2$,
(ii) $I=(G,H)+I^2$, and (iii) $I$ contains a monic, namely, $T-1$. Assume, if possible, that
there exist $g,h\in I$ such that $I=(g,h)$ with $g-G\in I^2$, and $h-H\in I^2$. Then, by 
\cite[2.14 (iv)]{brs1} there exists a matrix $\sigma\in SL_2(R[T])$ \emph{which is a lift of} 
$\ol{\sigma}$, and $(f,T-1)\sigma=(g,h)$. But $\sigma\in SL_2(R[T])=SL_2(k[X,Y,Z,T])$. As 
$SK_1(R[T])=0$, the image $\ol{\sigma}$  of $\sigma$ cannot be a non-trivial element of 
$SK_1(B)$. This is a contradiction. Therefore, the generators $\ol{G},\ol{H}$ of $I/I^2$ cannot be lifted to a
set of generators of $I$. 
\end{example}

\medskip

\remove{
\rmk
In the example given above note that $I(0)=R$ and  by \cite[3.9]{brs1}, $\ol{G},\ol{H}$ can actually be lifted 
to a set of generators of $I/(I^2T)$. Therefore (\ref{sri}) also shows that the results of Mandal \cite{m1, m2} and 
Bhatwadekar-Keshari \cite[4.13]{bk} cannot be improved in general for local complete intersection ideals of 
height $2$. 
}

\bc
Let $I$ be a non-zero ideal in $\k[X_1,\cdots, X_d]$ with $\mu(I/I^2)=n$. Then, $\mu(I)=n$  if $n\leq 2$ or $n\geq d-1$.
\ec

Let us now address the following question:

\smallskip

\quest\label{epi}
\emph{ Let $k$ be a field and $\phi:k[X_1,\cdots,X_n]\sur k[Y_1,\cdots,Y_m]$ be a surjective $k$-algebra morphism and let $I$ be the kernel of $\phi$.  Then  is $I$ generated by $n-m$ elements?} 

\smallskip

When $k$ is an algebraically closed field of characteristic zero, this is a famous open problem in Affine Algebraic Geometry (see \cite[Question 2.1]{dug} and the discussion surrounding it for an excellent survey).

For \emph{any} field $k$, there are  affirmative answers when: (1) $n=m+2$, or (2) $n\geq 2m+2$.  
Let us quickly indicate the solutions below. Before proceeding, note that under the hypothesis of
(\ref{epi}), $\mu(I/I^2)=n-m$.

\smallskip

\noindent
\emph{Proof of (1).} $n=m+2$. Then $\mu(I/I^2)=n-m=2$. Exactly the same argument as in (\ref{monic}) (Case 1), shows that
$I$ is surjective image of a projective $k[X_1,\cdots,X_n]$-module $P$ of rank $2$. But $P$ is free by the Quillen-Suslin Theorem. \qed

\smallskip

\noindent
\emph{Proof of (2).}   $n\geq 2m+2$. Note that $\mu(I/I^2)=n-m\geq m+2=\dd(k[X_1,\cdots,X_n]/I)+2$. Apply Mohan Kumar's result ((\ref{mohan1}) above). \qed

\medskip

When $k=\k$, we can do a little better, as the following theorem shows.

\bt
Let $\phi:\k[X_1,\cdots,X_n]\sur \k[Y_1,\cdots,Y_m]$ be a surjective $\k$-algebra morphism and let $I$ be the kernel of $\phi$. Assume that $n\geq 2m+1$. Then $\mu(I)=n-m$. 
\et

\proof
 $\mu(I/I^2)=n-m\geq m+1=\dd(\k[X_1,\cdots,X_n]/I)+1$. Apply (\ref{murthy}).
\qed

\medskip

\rmk
Let $k=\k$ and $m=2$. Then (\ref{epi}) has an affirmative answer for all
$n\geq 3$.

\section{Addition and subtraction principles}\label{AddSub}
Adapting the arguments from \cite[3.2, 3.3]{brs3}, in this section we prove the following 
\emph{``addition"} and \emph{``subtraction"} principles.
We shall need them in the next section. 

\bt\label{add}(Addition principle) Let $R$ be an affine algebra of dimension $d\geq 3$ over $\k$. Let 
$I,J$ be two comaximal ideals of $R$, each  of height $d-1$. Assume that $I=(a_1,\cdots,a_{d-1})$ and 
$J=(b_1,\cdots,b_{d-1})$. Then,
$I\cap J=(c_1,\cdots,c_{d-1})$ such that $c_i-a_i\in I^2$ and $c_i-b_i\in J^2$ for $i=1,\cdots,d-1$.
\et

\proof 
For $d\geq 5$, one may  appeal to \cite[3.1]{brs4}. The following proof works for $d\geq 3$.

Let `bar' denote reduction modulo $J$.
Note that $(\ol{a}_1,\cdots,\ol{a}_{d-1})$ is a unimodular row in $R/J$ and can be completed to a matrix $\tau\in SL_{d-1}(R/J)$. Since $SL_{d-1}(R)\lra SL_{d-1}(R/J)$ is surjective by (\ref{vas}), we can lift $\tau$ to $\gamma\in SL_{d-1}(R)$, replace $(a_1,\cdots,a_{d-1})$ by $\gamma(a_1,\cdots,
a_{d-1})$ and assume that $(a_1,\cdots,a_{d-2})+J=R$ and $a_{d-1}\in J$.
Adding  suitable multiples of $a_{d-1}$ to $a_1,\cdots,a_{d-2}$, we can further ensure that $\hh(a_1,\cdots,a_{d-2})=d-2$. 

 Let $K=(a_1,\cdots,a_{d-2})$ and let `tilde' denote reduction modulo  $K$. As $K+J=R$, the row
$(\widetilde{b}_1,\cdots,\widetilde{b}_{d-1})$ is unimodular over $R/K$. As $\dd(R/K)\leq 2$, the stably free 
$R/K$-module defined by this unimodular row is free by (\ref{free}). In other words, 
$(\widetilde{b}_1,\cdots,\widetilde{b}_{d-1})$ is completable over $R/K$. As a consequence, the unimodular row 
$(b_1,\cdots,b_{d-1})_{1+K}$ over the ring $R_{1+K}$ is completable to a matrix in 
$SL_{d-1}(R_{1+K})$. We can clear denominators and find
a suitable $s$ of the form $1+t$, where $t\in K$, and some $\alpha\in SL_{d-1}(R_{s})$ such that 
$(b_1,\cdots,b_{d-1})_{s}\alpha=(1,0,\cdots,0)$.
Clearly, we can adjust so that $s\in J$.

We have $(I\cap J)_s=I_s=(a_1,\cdots,a_{d-1})_s$. Also, 
$(I\cap J)_t=J_t=(b_1,\cdots,b_{d-1})_t$.

As $t\in K=(a_1,\cdots,a_{d-2})$,  the unimodular row  
$(a_1,\cdots,a_{d-1})_t$ contains a subrow of shorter length, namely, $(a_1,\cdots,a_{d-2})_t$.
Therefore, $(a_1,\cdots,a_{d-1})_t$ is elementarily completable. In other words, there is 
$\sigma\in E_{d-1}(R_t)$ such that 
$(a_1,\cdots,a_{d-1})_t\sigma=(1,0,\cdots,0)$.

Over the ring $R_{st}$ we have
$(a_1,\cdots,a_{d-1})\sigma_s\alpha_t^{-1}=(b_1,\cdots,b_{d-1})$ and therefore, 
$$(a_1,\cdots,a_{d-1})\alpha_t^{-1}(\alpha_t\sigma_s\alpha_t^{-1})=(b_1,\cdots,b_{d-1}).$$
Write $\theta=\alpha_t\sigma_s\alpha_t^{-1}$. Note that $\sigma_s$ is isotopic to identity and
consequently so is $\theta$. By  (\ref{qsplit}), 
$\theta=\theta'_t\theta''_s$, where $\theta'\in SL_{d-1}(R_s)$ and $\theta''\in SL_{d-1}(R_t)$.
Let $(a_1,\cdots,a_{d-1})\alpha^{-1}\theta'=(x_1,\cdots,x_{d-1})$ and 
$(b_1,\cdots,b_{d-1}){\theta''}^{-1}=(y_1,\cdots,y_{d-1})$.
Then we have $(I\cap J)_s=(x_1,\cdots,x_{d-1})$ and $(I\cap J)_t= (y_1,\cdots,y_{d-1})$. As these generators agree over 
$R_{st}$, we can patch them to obtain a set of generators of $I\cap J$, say, $I\cap J=(z_1,\cdots,z_{d-1})$.

Recall that  `tilde' denotes reduction modulo $I$ and `bar' denotes that modulo $J$.
Observe that by the above construction, $(\wt{a}_1,\cdots,\wt{a}_{d-1})$ and $(\wt{z}_1,\cdots,
\wt{z}_{d-1})$ differ by an 
element of $SL_{d-1}(R/I)$. Similarly, $(\ol{b}_1,\cdots,\ol{b}_{d-1})$ and $(\ol{z}_1,\cdots,\ol{z}_{d-1})$
differ by an element of $SL_{d-1}(R/J)$. Now we use $SL_{d-1}(R)\sur SL_{d-1}(R/I\cap J)\iso SL_{d-1}(R/I)\times SL_{d-1}(R/J)$
to get a suitable $\beta\in SL_{d-1}(R)$ and observe that if $(z_1,\cdots,z_{d-1})\beta=(c_1,\cdots,c_{d-1})$, then these 
are the generators of $I\cap J$ we were looking for.
\qed

\bt\label{sub2}(Subtraction principle) Let $R$ be an affine algebra of dimension $d\geq 3$ over 
$\k$ (with the additional assumption $p\neq 2$ if $d=4$). Let $I,J$ be two comaximal ideals of $R$ such that $\hh(I)=d-1$ and $\hh(J)\geq d-1$. Let $P$ be a projective $R$-module of rank $d-1$ with trivial determinant and $\chi:R\iso \wedge^{d-1}P$ be an isomorphism. Let 
$\ga:P\sur I\cap J$ and $\gb:R^{d-1}\sur I$ be surjections. Suppose that there exists an isomorphism
$\delta:(R/I)^{d-1}\iso P/IP$ such that: 
\begin{enumerate}
\item
$\wedge^{d-1}\delta=\chi\ot R/I$,
\item
$(\ga\ot R/I)\delta=\gb\ot R/I$.
\end{enumerate} 
Then, there exists a surjection $\gamma:P\sur J$ such that $\gamma\ot R/J=\ga\ot R/J$.
\et

\proof
We first remark that since $J$ is locally generated by $d-1$ elements, either it is proper of height 
$d-1$, or $J=R$.

Let $\gb$ correspond to $I=(a_1,\cdots,a_{d-1})$. As before, we first note that we can always make changes up to $SL_{d-1}(R)$ in our arguments. We can make similar reductions as in the above proof and assume that: (1) $(a_1,\cdots,a_{d-2})+J^2=R$, 
(2) $a_{d-1}\in J^2$, and (3) $\hh(a_1,\cdots,a_{d-2})=d-2$. Once these are done, we pick some
$\lambda\in (a_1,\cdots,a_{d-2})$ such that $\lambda\equiv 1$ modulo $J^2$, replace $a_{d-1}$ by
$\lambda+a_{d-1}$ and obtain $a_{d-1}\equiv 1$ modulo $J^2$.

Consider the following ideals in $R[T]$:  $K'=(a_1,\cdots,a_{d-2},T+a_{d-1})$, $K''=J[T]$, and $K=K'\cap K''$. We aim to show that there is a surjection $\theta:P[T]\sur K$ such that $\theta(0)=\ga$. If we can achieve so then we
can specialize at $1-a_{d-1}$ to obtain $\gamma:=\theta(1-a_{d-1}):P\sur J$. As $a_{d-1}\equiv 1$ modulo $J^2$, we have $\gamma\ot R/J=\theta(1-a_{d-1})\ot R/J=\theta(0)\ot R/J=\ga\ot R/J$, and we will be done. Rest of the proof is about finding such a $\theta$, and we break it into two cases.

\smallskip

\noindent
\emph{Case 1.} Assume that $d\geq 5$. As $\dd(R[T]/K')\leq 2$, the module $P[T]/K'P[T]$ is free of rank $d-1$ by (\ref{free}). We choose an isomorphism $\kappa(T):(R[T]/K')^{d-1}\iso P[T]/K'P[T]$ such that
$\wedge^{d-1}\kappa(T)=\chi\ot R[T]/K'$. Therefore, $\wedge^{d-1}\kappa(0)=\chi\ot R/I=\wedge^{d-1}\delta$ and we conclude that $\kappa(0)$ and $\delta$ differ by an element $\mathfrak{a}\in SL_{d-1}(R/I)$. By (\ref{vas}), we can find a lift of $\mathfrak{a}$ in $SL_{d-1}(R)$ and use this to alter $\kappa(T)$ so that $\kappa(0)=\delta$.

Sending the canonical basis vectors to $a_1,\cdots,a_{d-2},T+a_{d-1}$ respectively, we have a surjection
from $R[T]^{d-1}$ to $K'$, which will induce a surjection $\epsilon(T):(R[T]/K')^{d-1}\sur K'/K'^2$.
Composing with $\kappa(T)^{-1}$ we have:
$$\pi(T):=\epsilon(T)\kappa(T)^{-1}:P[T]/K'P[T]\sur K'/K'^2.$$
Observe that $\pi(0)=\epsilon(0)\kappa(0)^{-1}=(\beta\ot R/I)\delta^{-1}=\ga\ot R/I$ (see hypothesis
(2) above). We can now apply \cite[2.3]{mrs} to obtain a surjection $\theta(T):P[T]\sur K$ so that
$\theta(0)=\ga$.

\smallskip

\noindent
\emph{Case 2.} Assume that $d=3\text{ or } 4$.
Write $L=(a_1,\cdots,a_{d-2})$.
 Note that $\dd(R/L)=2$,  and therefore, $P/LP$ is free by (\ref{free}). Consequently, $P_{1+L}$ is a free 
$R_{1+L}$-module of rank $d-1$. We choose an isomorphism $\xi:R_{1+L}^{d-1}\iso P_{1+L}$ 
such that $\wedge^{d-1}\xi=\chi\ot R_{1+L}$. Note that $\xi$ induces isomorphism 
$\xi(T): R_{1+L}[T]^{d-1}\iso P_{1+L}[T]$.

We have $K_{1+L}=K'_{1+L}$.
Sending the canonical  basis vectors to $a_1,a_2,\cdots, a_{d-2},T+a_{d-1}$, respectively,  we have a surjection 
$\pi(T):R_{1+L}[T]^{d-1}\sur K_{1+L}$.
We therefore have
$$\gamma(T):=\pi(T)(\xi(T))^{-1}:P_{1+L}[T]\sur K_{1+L}.$$
Recall that we have surjections $\ga:P\sur I\cap J$ and $\gb:R^{d-1}\sur I$. We would like to compare 
$\gamma(0)$ 
and $\ga_{1+L}:P_{1+L}\sur I_{1+L}$.

Note that $\gamma(0)=\pi(0)\xi^{-1}=\gb_{1+L}\xi^{-1}$ and $I_{1+L}/I^2_{1+L}=I/I^2$. 
We have induced surjections:
$$\ol{\gb}\ol{\xi}^{-1}:P/IP\iso (R/I)^{d-1}\sur I/I^2$$
$$\ol{\gb}\gd^{-1}:P/IP\iso (R/I)^{d-1}\sur I/I^2$$

We have,  $\wedge^{d-1} \gd=\chi\ot R/I$ and $\wedge^{d-1} \xi=\chi\ot R_{1+L}$. 
Further, 

\begin{enumerate}
\item
$\gamma(0)\ot R/I=(\gb\ot R/I)(\xi^{-1}\ot R/I)$,
\item
$\ga_{1+L}\ot R/I=(\gb\ot R/I)\gd^{-1}$. 
\end{enumerate}
It follows that $\gamma(0)\ot R/I$ and $\ga_{1+L}\ot R/I$ differ by an element  
$\mathfrak{a}\in SL(P/IP)$. As $P/IP$ is free and $\dd(R/I)=1$, by (\ref{vas}), 
$SL(P/IP)=E(P/IP)$ for $d= 4$.
As $E(P_{1+L})\lra E(P_{1+L}/(IP)_{1+L})=E(P/IP)$ is surjective, we can find a lift of
$\mathfrak{a}$ in $E(P_{1+L})$. On the other hand, if $d=3$, 
applying \cite[2.3]{brs3} we conclude that $\gamma(0)$ and $\ga_{1+L}$ differ by an automorphism in 
$SL(P_{1+L})$.  In either case, we can alter $\gamma(T)$ and assume that $\gamma(0)=\ga_{1+L}$. 
We can find
a suitable $t\in L$ such that: (a) $1+t\in J$, 
 (b) $P_{1+tR}$ is free, and (c) there is a surjection 
$\gamma(T):P_{1+tR}[T]\sur I_{1+tR}$ such that $\gamma(0)=\ga_{1+tR}$.

On the other hand, $\ga_t:P_t\sur J_t$ will induce 
$\ga_{t}(T):P_{t}[T]\sur K_{t}$ (note that $K_t=K''_t=J_t[T]$). Now, the two surjections
$$\gamma(T)_{t}:P_{t(1+tR)}[T]\sur K_{t(1+tR)}[T]=R_{t(1+tR)}[T]$$ 
$$\ga_{t(1+tR)}:P_{t(1+tR)}[T]\sur K_{t(1+tR)}=R_{t(1+tR)}[T]$$ 
are such that they agree when $T=0$, and  both their kernels are:
\begin{enumerate}
\item
free, if $d=3$ (as any unimodular row of length $2$ is completable),
\item
extended, if $d=4$ and $2R=R$ (by (\ref{ravi})). 
\end{enumerate}

A standard patching argument using (\ref{qsplit}) yields a surjection $\theta:P[T]\sur K$
such that $\theta(0)=\ga$. 
\qed

\medskip

Taking $P=R^{d-1}$ in (\ref{sub2}), we get the following corollary.
  
\bc\label{sub}(Subtraction principle) 
Let $R$ be an affine algebra of dimension $d\geq 3$ over 
$\k$ (with the additional assumption $p\neq 2$ if $d=4$). Let $I,J$ be two comaximal ideals of 
$R$, each  of height $d-1$. Assume that $I=(a_1,\cdots,a_{d-1})$ and $I\cap J=(c_1,\cdots,c_{d-1})$ such that $c_i-a_i\in I^2$
for $i=1,\cdots,d-1$. 
Then, $J=(b_1,\cdots,b_{d-1})$ such that   $c_i-b_i\in J^2$ for $i=1,\cdots,d-1$.
\ec

\smallskip

Taking $J=R$ in (\ref{sub2}), we get the following corollary.

\bc\label{sub3}
Let $R$ be an affine algebra of dimension $d\geq 3$ over 
$\k$ (with the additional assumption $p\neq 2$ if $d=4$). Let $I$ be an ideal of $R$ such that $\hh(I)=d-1$. Let $P$ be a projective $R$-module of rank $d-1$ with trivial determinant and $\chi:R\iso \wedge^{d-1}P$ be an isomorphism. Let 
$\ga:P\sur I$ and $\gb:R^{d-1}\sur I$ be surjections. Suppose that there exists an isomorphism
$\delta:(R/I)^{d-1}\iso P/IP$ such that: 
\begin{enumerate}
\item
$\wedge^{d-1}\delta=\chi\ot R/I$,
\item
$(\ga\ot R/I)\delta=\gb\ot R/I$.
\end{enumerate} 
Then, there exists a surjection $\gamma:P\sur R$.
\ec

\smallskip

\rmk
When $p\neq 2$ and $R$ is smooth,  later in (\ref{addsub}) we prove: If $I,J\subset R$ are 
comaximal ideals of height $d-1$, and if two of $I$, $J$, $I\cap J$ are complete intersections,
then so is the third.

\remove{
We end this section with a result which is easy to prove using (\ref{move}), (\ref{add}), and 
(\ref{sub}). To get an idea of the proof, see \cite[3.3]{brs4} or \cite[4.1]{brs3}.

\bp\label{move1}
Let $R$ be an affine algebra of dimension $d\geq 3$ over 
$\k$, with the following additional assumptions: (i) $2R=R$ if $d=4$, (ii) $R$ is smooth if $d\geq 5$.  Let $J,J_1,J_2$ be ideals of $R$, each of height $d-1$, such that $J+J_1\cap J_2=R$. Assume that 
$J\cap J_1=(a_1,\cdots,a_{d-1})$ and $J\cap J_2=(b_1,\cdots,b_{d-1})$, where $a_i-b_i\in J^2$ for 
$i=1,\cdots,d-1$. Suppose that there is an ideal $J_3$ of height $d-1$ such that $J_3+J\cap J_1\cap J_2=R$ and $J_3\cap J_1=(c_1,c_2)$
with $c_i-a_i\in J_1^2$ for $i=1,\cdots,d-1$. Then  $J_3\cap J_2=(d_1,\cdots,d_{d-1})$
with $d_i-b_i\in J_2^2$ for $i=1,\cdots,d-1$. 
\ep
}

\section{An obstruction group and an obstruction class}\label{ECG}
Let $R$ be an affine algebra of dimension $d\geq 3$ over $\k$.
Let $I\subset R$ be an ideal of height $d-1$ such that $I/I^2$ is generated by $d-1$ elements. Two
surjections $\ga,\gb:(R/I)^{d-1}\sur I/I^2$ are said to be \emph{related} if there exists 
$\sigma\in SL_{d-1}(R/I)$ such that $\ga\sigma=\gb$. This defines an equivalence relation on 
the set of surjections from $(R/I)^{d-1}$ to $I/I^2$. 

\bl\label{equi}
Let $R$ be an affine algebra of dimension $d\geq 3$ over $\k$.
Let $I\subset R$ be an ideal of height $d-1$ such that $\mu(I/I^2)=d-1$.
Let a surjection $\ga:(R/I)^{d-1}\sur I/I^2$ be such that it has a surjective lift 
$\theta : R^{d-1}\sur I$.
Then the same is true for any $\gb$ related to $\ga$.
\el

\proof
Let $\ol{\sigma}\in SL_{d-1}(R/I)$ be such that $\ga\ol{\sigma}=\gb$. By the hypothesis, there is a
surjection $\theta:R^{d-1}\sur I$ such that $\theta\ot R/I=\ga$.
 Applying Theorem \ref{vas} we can  find a lift $\sigma\in SL_{d-1}(R)$ of $\ol{\sigma}$.
Then $\theta\sigma: R^{d-1}\sur I$ is a surjective lift of $\gb$.
\qed

\bd
An equivalence class of surjections from $(R/I)^{d-1}$ to $I/I^2$ will be called a \emph{local orientation} of $I$. A local orientation $\omega_I$ of $I$ will be called a \emph{global orientation} 
if  a surjection (hence all) in the class of $\omega_I$  can be lifted to a surjection from $R^{d-1}$ to $I$.
\ed

\medskip

\rmk
Let $R$ be an affine algebra of dimension $d\geq 3$ over $\k$.
We now define the $(d-1)$\textsuperscript{th} Euler class group of $R$. Recall that for any commutative Noetherian ring $A$ of dimension $d$ and any integer $n$ with $2n\geq d+3$, there
is a notion of the $n$\textsuperscript{th} Euler class group $E^n(A)$ of $A$ in \cite{brs4}.
Our definition is modeled on their definition but there is a difference:  we consider 
the action of $SL_{d-1}(R/I)$ as defined above, whereas they consider the action of $E_{d-1}(R/I)$. This difference is crucial.

\medskip

\bd\label{def}
\emph{\underline{The Euler class group $E^{d-1}(R)$}:} Let $R$ be an affine algebra of dimension $d\geq 3$ over $\k$.
Let $G$ be the free abelian group on the set $B$ of pairs $({\mathcal I},\omega_{\mathcal I})$ where 
${\mathcal I}\subset R$ is an ideal of
height $d-1$ with the property that $\text{Spec}\,(R/{\mathcal I})$ is connected,   
$\mu({\mathcal I}/{\mathcal I}^2)=d-1$,
and $\omega_{\mathcal I}:(R/{\mathcal I})^{d-1}\sur {\mathcal I}/{\mathcal I}^2$ is a local orientation.

Let $I$ be any ideal of $R$ of height $d-1$ such that $I/I^2$ is generated by $d-1$ elements. Then $I$ has a unique decomposition, $I={\mathcal I}_1\cap\cdots\cap {\mathcal I}_k$, where each of 
$\text{Spec}(R/{\mathcal I}_i)$ is connected, ${\mathcal I}_i$ are pairwise comaximal, and 
$\text{ht}\,{\mathcal I}_i=d-1$ (for a
proof see \cite{brs3} or \cite{d}). Now if $\omega_I$ is a local orientation of $I$,
then it naturally gives rise to $\omega_{{\mathcal I}_i}:(R/{\mathcal I}_i)^{d-1}\sur {\mathcal I}_i/{{\mathcal I}_i}^2$ for $1\leq i\leq k$.  
By $(I,\omega_I)$ we mean the element $\sum ({\mathcal I}_i,\omega_{{\mathcal I}_i})\in G$.

Let $H$ be the subgroup of $G$ generated by the set $S$ of pairs $(I,\omega_I)$ in $G$ such that 
$\omega_I$ is a global orientation. We define the $(d-1)$\textsuperscript{th} Euler class group of 
$R$ as $E^{d-1}(R):=G/H$.  
\ed

\medskip

\rmk
In view of (\ref{equi}), we would not distinguish between a surjection $\omega_I:(R/I)^{d-1}\sur I/I^2$
and the equivalence class it represents.

We now state a lemma from \cite[4.1]{ke}.

\begin{lemma}\label{mphil}
Let $G$ be a free abelian group with basis $B=(e_i)_{i\in {\mathcal I}}$ . Let $\sim$  be an equivalence relation on
$B$. Define $x\in G$  to be ``reduced" if $x=e_1+\cdots +e_r$ and $e_i\not=e_j$ for $i \not=j$.
Define $x\in G$ to be ``nicely reduced" if $x=e_1+\cdots +e_r$ is such that  $e_i\not\sim e_j$ for $i \not=j$.
Let $S\subset G$ be such that 
\begin{enumerate}
\item
Every element of $S$ is nicely reduced.
\item
Let $x,y\in G$ be such that each of $x,y,x+y$ is nicely reduced. If two of $x,y,x+y$ are in $S$, then so is the third.
\item
Let $x\in G\setminus S$ be nicely reduced and let ${\mathcal J}\subset {\mathcal I}$ be  finite. Then there exists $y\in G$ with the following properties : (i) $y$ is nicely reduced; (ii) $x+y\in S$; (iii) $y+e_j$ is nicely reduced $\forall j\in{\mathcal J}$. 
\end{enumerate}

Let $H$ be the subgroup of $G$ generated by $S$. If $x \in H$ is nicely reduced, then $x \in  S$.
\end{lemma}

We are now ready to prove:

\bt\label{zero}
Let $R$ be an affine algebra of dimension $d\geq 3$ over 
$\k$ (with the additional assumption $p\neq 2$ if $d=4$).
Let $I$ be any ideal of $R$ of height $d-1$ such that $\mu(I/I^2)=d-1$, and 
$\omega_I:(R/I)^{d-1}\sur I/I^2$ be a local orientation.
Assume that the image of $(I,\omega_I)$ is trivial in $E^{d-1}(R)$. Then $\omega_I$ is global.
In other words, there is an
$R$-linear  surjection $\theta : (R/I)^{d-1}\sur I$ such that $\theta $ lifts $\omega_I$.
\et

\proof
We take $G$ to be the free abelian group generated by $B$, as defined in (\ref{def}). Define a relation $\sim$ on $B$ as: $(K,\omega_K)\sim (K',\omega_{K'})$ if 
$K=K'$. Then it is an equivalence relation. 

Let $S\subset G$ be as in (\ref{def}). In view of the above lemma, it is enough to show that the three
conditions in (\ref{mphil}) are satisfied. Condition (1) is clear, almost from the definition.
The addition and subtraction principles (\ref{add}, \ref{sub}) will yield  condition (2). Finally, applying the moving lemma (\ref{move}), it is clear that (3) is also satisfied.
\qed

\medskip

\notation
 Let 
$\sigma\in GL_{d-1}(R/I)$ and let $\text{det}(\sigma)=\ol{u}$. As we are dealing with the action
of $SL_{d-1}(R/I)$ on  surjections from $(R/I)^{d-1}$ to $I/I^2$ in the definitions above, there will be no ambiguity if
we write $\omega\sigma$ as $\ol{u}\omega$ and the corresponding element in $E^{d-1}(R)$ as $(I,\ol{u}\omega)$. 

\medskip

\rmk\label{smb}
Let $\omega$, $\omega'$ be two surjections from $(R/I)^{d-1}$ to $I/I^2$. 
With the above notations in mind,  $(I,\omega')=(I,\ol{u}\omega)$ for some unit $\ol{u}$ in $R/I$
by \cite[2.2]{b}.

\smallskip

We shall need the following result in Section \ref{splitting}.

\bp\label{square}
Let $R$ be as in (\ref{def}) and $(I,\omega)\in E^{d-1}(R)$. Let $\ol{u}$ in $R/I$ be a unit.
The following assertions hold:
\begin{enumerate}
\item
$(I,\omega)=0$ if and only if $(I,\ol{u}^2\omega)=0$ in $E^{d-1}(R)$.
\item
$(I,\omega)=(I,\ol{u}^2\omega)$ in $E^{d-1}(R)$.
\end{enumerate}
\ep

\proof
Statement (1) is needed to prove (2). The proofs are the same as in \cite[5.3, 5.4]{brs3} and we do not repeat them. One has to simply apply (\ref{vas}) at the appropriate places to lift automorphisms.
\qed

\bd\label{eulerclass}
\emph{\underline{The $(d-1)$\textsuperscript{th} Euler class of a projective module}:} Let $R$ be an affine algebra of dimension $d\geq 3$ over 
$\k$, with the following additional assumptions: (i) $p\neq 2$ if $d=4$, (ii) $R$ is smooth if $d\geq 5$. Let $P$ be a projective $R$-module of rank $d-1$ with trivial determinant. Fix an isomorphism $\chi: R\iso \wedge^{d-1}(P)$. Let $\alpha:P\sur J$ be a surjection where $J\subset R$ is an ideal of height $d-1$. Note that $P/JP$ is a free $R/J$-module.
Choose an isomorphism $\sigma:(R/J)^{d-1}\iso P/JP$ such that $\wedge^{d-1}\sigma=\chi\ot R/J$. Let $\omega_J$ be the composite:
$$(R/J)^{d-1}\stackrel{\sigma}{\iso} P/JP\stackrel{\ol{\ga}}{\sur} J/J^2$$
We define the \emph{$(d-1)$\textsuperscript{th} Euler class} of the pair $(P,\chi)$ as $e_{d-1}(P,\chi):=(J,\omega_J)$ (we show below that this association does not depend on the chosen surjection $\alpha: P\sur J$). We shall simply call it the Euler class here.
\ed

\bt
Let $R$ be an affine algebra of dimension $d\geq 3$ over 
$\k$, with the following additional assumptions: (i) $p\neq 2$ if $d=4$, (ii) $R$ is smooth if $d\geq 5$. Let $P$ be a projective $R$-module of rank $d-1$ with trivial determinant. Fix an isomorphism $\chi:R\iso \wedge^{d-1} P$.
The Euler class $e_{d-1}(P,\chi)$ is well-defined.
\et

\proof
Let $\beta:P\sur J'$ be another surjection such that $J'$ is an ideal of $R$ of height $d-1$.

By (\ref{hom}) there exists an ideal $I\subset R[T]$ of height $d-1$ and a surjection $\phi(T):P[T]\sur I$ such that $I(0)=J$,
$\phi(0)=\ga$ and $I(1)=J'$,  $\phi(1)=\gb$. 

Let $N=(I\cap R)^2$. Then $\hh(N)\geq d-2$ and therefore, $\dd(R/N)\leq 2$. By (\ref{free}), $P/NP$ is a free $R/N$-module of rank $d-1$. On the other hand, by the same reasoning, $P[T]/IP[T]$ is a free $R[T]/I$-module of rank $d-1$.

We can choose an isomorphism $\tau:(R/N)^{d-1}\iso P/NP$ such that $\wedge^{d-1}\tau=\chi\otimes R/N$.
This choice of $\tau$ gives us a basis of $P/NP$, which in turn induces  a basis of the free module 
$P[T]/IP[T]$. Using this basis of $P[T]/IP[T]$ and the surjection $\phi(T):P[T]\sur I$, we obtain a surjection $\omega:(R[T]/I)^{d-1}\sur I/I^2$. Note that, due to the choice of the basis, $\omega(0)=\omega_J:(R/J)^{d-1}\sur J/J^2$, and $\omega(1)=\omega_{J'}:(R/J')^{d-1}\sur J'/J'^2$.

Now using (\ref{chi}), we obtain an ideal $K$ of height $d-1$ with $K+J=R=K+J'$, and a 
surjection $\omega_K:(R/K)^{d-1}\sur K/K^2$ such that
$$(J,\omega_J)+(K,\omega_K)=(J',\omega_{J'})+(K,\omega_K) \text{ in } E^{d-1}(R).$$
This proves that $e_{d-1}(P,\chi)$ is well-defined.
\qed

\smallskip

We now prove that $e_{d-1}(P,\chi)$ is the precise obstruction for $P$ to split off a free summand of rank one.

\bt\label{split}
Let $R$ be an affine algebra of dimension $d\geq 3$ over 
$\k$, with the following additional assumptions: (i) $p\neq 2$ if $d=4$, (ii) $R$ is smooth if $d\geq 5$. Let $P$ be a projective $R$-module of rank $d-1$ with trivial determinant. Fix an isomorphism $\chi:R\iso \wedge^{d-1} P$. Then $P\iso Q\oplus R$ for some $R$-module $Q$ if and only if $e_{d-1}(P,\chi)=0$ in $E^{d-1}(R)$.
\et

\proof
Let $\alpha:P\sur J$ be a surjection such that $J\subset R$ is an ideal of height $d-1$. 
Choose an isomorphism $\sigma:(R/J)^{d-1}\iso P/JP$ such that $\wedge^{d-1}\sigma=\chi\ot R/J$. Let $\omega_J$ be the composite:
$$(R/J)^{d-1}\stackrel{\sigma}{\iso} P/JP\stackrel{\ol{\ga}}{\sur} J/J^2$$
Then, by definition (\ref{eulerclass}), $e_{d-1}(P,\chi)=(J,\omega_J)$ in $E^{d-1}(R)$.

Assume first that  $P\iso Q\oplus R$ for some $R$-module $Q$. Write $\alpha=(\beta,a)$. Applying an elementary automorphism of $P$, we may assume that height of $K:=\beta(Q)$ is at least $d-2$. Note that
$J=(K,a)$.

As the determinant of $Q$ is also trivial, we may assume that $\chi$ is induced by 
$\chi':R\iso \wedge^{d-2} Q$. As $\dd (R/K)\leq 2$, by (\ref{free}) the $R/K$-module $Q/KQ$ is free of rank $d-2$.
Choose an isomorphism $\gamma':(R/K)^{d-2}\iso Q/KQ$ such that $\wedge^{d-2} \gamma'=\chi'\ot R/K$.
We have a surjection $(\beta\ot R/K)\gamma': (R/K)^{d-2}\sur K/K^2$. Suppose that this induces 
$K=(\ga_1,\cdots,\ga_{d-2})+K^2$. Then, by (\ref{nak}) there exists $e\in K^2$ such that $K=(\ga_1,\cdots,\ga_{d-2},e)$ where $e(1-e)\in (\ga_1,\cdots,\ga_{d-2})$. As $J=(K,a)$, it follows from 
(\ref{nak} (3)) that $J=(\ga_1,\cdots,\ga_{d-2},b)$, where $b=e+(1-e)a$. Note that 
$b-a=e-ea\in K^2\subset J^2$. It is now easy to see that $\omega_J$ is induced by $\ga_1,\cdots,\ga_{d-2},b$.  In other words, $(J,\omega_J)=0$ in $E^{d-1}(R)$.

Conversely, assume that $e_{d-1}(P,\chi)=0$. Then $(J,\omega_J)=0$ in $E^{d-1}(R)$ and therefore by (\ref{zero}) there exists $\theta:R^{d-1}\sur J$ which is a lift of $\omega_J$. Let $\omega_J$ correspond to 
$J=(a_1,\cdots,a_{d-1})+J^2$ and $\theta$ correspond to $J=(b_1,\cdots,b_{d-1})$. We have
$b_i-a_i\in J^2$ for $1\leq i\leq d-1$. 
Now apply (\ref{sub3}).
\qed

\bt\label{ind}
Let $R$ be an affine algebra of dimension $d\geq 3$ over 
$\k$, with the following additional assumptions: (i) $p\neq 2$ if $d=4$, (ii) $R$ is smooth if $d\geq 5$. Let $P$ be a projective $R$-module of rank $d-1$ with trivial determinant. Fix an isomorphism 
$\chi:R\iso \wedge^{d-1} P$. Let $e_{d-1}(P,\chi)=(J,\omega_J)$ for some $(J,\omega_J)\in E^{d-1}(R)$. Then there is a surjection $\alpha:P\sur J$ such that $(J,\omega_J)$ is induced by $(\alpha,\chi)$.
\et

\proof
We choose an isomorphism $\sigma:(R/J)^{d-1}\iso P/JP$ such that $\wedge^{d-1}\sigma=\chi\ot R/J$. 
We have the composite:
$$P/JP\stackrel{\sigma^{-1}}{\iso} (R/J)^{d-1}\stackrel{\omega_J}{\sur} J/J^2,$$
let us call it $\ol{\theta}$. Applying (\ref{move}), we can find an ideal $J'\subset R$ and a surjection
$\eta:P\sur J\cap J'$ such that: (i) $J+J'=R$, (ii) $\hh(J')\geq d-1$, and
(iii) $\eta\ot R/J=\ol{\theta}$. If $J'=R$, we are done. Let $\hh(J')=d-1$. It is easy to see that the pair 
$(\eta,\chi)$ induces $e_{d-1}(P,\chi)=(J,\omega_J)+(J',\omega_{J'})$ in $E^{d-1}(R)$.
As $e_{d-1}(P,\chi)=(J,\omega_J)$, we have $(J',\omega_{J'})=0$ in $E^{d-1}(R)$.
By (\ref{zero}), there is a surjection $\beta:R^{d-1}\sur J'$ which is a lift of 
$\omega_{J'}$.  Now apply (\ref{sub2}).
\qed

\smallskip

\rmk
We reiterate that for the definition and the results involving only $E^{d-1}(R)$, we do not need any smoothness assumption. We need smoothness whenever we talk about the Euler class of a projective module and $d\geq 5$. We could have imposed a blanket assumption that $R$ is smooth throughout this section (or the article) but we decided against it as some subtle (and perhaps useful) points will be lost.

\section{Applications I: splitting problem vis-\`{a}-vis complete intersections}\label{splitting}

In the preceding section we have established that the Euler class $e_{d-1}(P,\chi)$ is the precise obstruction for $P$ to split off a free summand of rank one. In this section we refine the conclusion further and  give it a much more simplified and tangible form.

Let $R$ be an affine algebra of dimension $d$ over $\k$ and $P$ be a projective $R$-module of rank $d-1$ with trivial determinant. Without getting into the Euler class theory, we (re)prove the following basic result, whose idea  is essentially from \cite[Theorem 1]{mo3}.

\bt\label{mohan}
Let $R$ be an affine algebra of dimension $d\geq 3$ over $\k$ and $P$ be a projective $R$-module of rank $d-1$ with trivial determinant. Let $\alpha:P\sur J$ be an $R$-linear surjection such that $\hh(J)=d-1$. Assume that $P\iso Q\oplus R$ for some $R$-module $Q$. Then $\mu(J)=d-1$.
\et

\proof
Note that the case $d=3$ is trivial. Therefore, assume that $d\geq 4$.
Let $\alpha_{|_Q}=\beta$, and $\ga(0,1)=a$. Using a standard general position argument we may assume that $\gb(Q)=K$ is such that $\dd(R/K)\leq 2$. 

We have the induced surjection $Q/KQ\sur K/K^2$.
As the determinant of $Q$ is trivial, by (\ref{free})  $Q/KQ$ is free of rank $d-2$. 
Therefore, $\mu(K/K^2)=d-2$, and  by (\ref{nak} (3)) we have $\mu(J)=d-1$.
\qed

\medskip

In this context, we can naturally ask the following question.

\smallskip

\quest\label{ci}
\emph{ Let $R$ be an affine algebra of dimension $d\geq 3$ over $\k$ and $P$ be a projective $R$-module of rank $d-1$ with trivial determinant. Let $\alpha:P\sur J$ be an $R$-linear surjection such that $\hh(J)=d-1$. Assume that $\mu(J)=d-1$. Then, is $P\iso Q\oplus R$ for some $R$-module $Q$?}

To answer this question we first need the following result.

\bt\label{cilift}
Let $p\neq 2$ and $R$ be a smooth affine algebra of dimension $d\geq 3$ over $\k$. Let $J\subset R$ be an ideal of height $d-1$ such that $\mu(J)=d-1$. Suppose, it is given that $J=(b_{1},\cdots,b_{d-1})+J^2$. Then,  there exist $c_{1},\cdots,c_{d-1}\in J$
such that $J=(c_{1},\cdots,c_{d-1})$ and $b_{i}\equiv c_{i}$ mod $J^2$ for $i=1,\cdots,d-1$.
\et

\proof
Let $J=(a_{1},\cdots,a_{d-1})$.
By \cite[2.2]{b}, there is a matrix $\delta\in GL_{d-1}(R/J)$ such that
$(\overline{a_{1}},\cdots,\overline{a_{d-1}})\delta =
(\overline{b_{1}},\cdots,\overline{b_{d-1}})$. Let $u\in R$ be such that 
$\overline{u}=$det\,$(\delta)^{-1}$. 
Using Swan's Bertini Theorem \cite[2.11]{brs2}, we 
may assume that $B=R/(a_3,\cdots,a_{d-1})$ 
is a smooth affine three-fold over $\k$. 
By \cite[7.5]{frs}, the unimodular row $(\widetilde{u},\widetilde{a_{1}},\widetilde{a_{2}})\in Um_{3}(B)$ is completable.
Applying \cite[2.4]{rs} we have a set of generators of $J$, say 
$J=(\beta_{1},\cdots,\beta_{d-1})$,
and a matrix $\ol{\theta}\in SL_{d-1}(R/J)$ such that
$(\overline{\beta_{1}},\cdots,\overline{\beta_{d-1}})\ol{\theta}=
(\overline{b_{1}},\cdots,\overline{b_{d-1}})$. We can now apply 
(\ref{vas}) and  lift $\ol{\theta} $
to a matrix $\theta\in SL_{d-1}(R)$. Suppose that
$(\beta_{1},\cdots,\beta_{d-1})\theta=(c_{1},\cdots,c_{d-1})$. Then 
$J=(c_{1},\cdots,c_{d-1})$, as  desired.
\qed

\medskip

The addition and subtraction principles (\ref{add}, \ref{sub}) proved earlier now take the following form.
 
\bc\label{addsub}
Let $p\neq 2$, and let $R$ be a smooth affine algebra of dimension $d\geq 3$ over $\k$. Let $I$,
$J$ be two comaximal ideals of $R$, each of height $d-1$. If two of $I,J,I\cap J$ are generated by $d-1$ elements, then so is the third. 
\ec

\medskip

The following corollary will be used in the next section.

\bc\label{equal}
Let $p\neq 2$, and let $R$ be a smooth affine algebra of dimension $d\geq 3$ over $\k$. Let 
$J\subset R$ be 
an ideal of height $d-1$ such that $\mu(J/J^2)=2$. Let $\omega_1,\omega_2$ be two surjections from 
$(R/J)^{d-1}$ to $J/J^2$. Then $(J,\omega_1)=(J,\omega_2)$ in $E^{d-1}(R)$.
\ec

\proof
If $(J,\omega_1)=0$, then the conclusion follows from the above theorem. Assume that 
$(J,\omega_1)\neq 0$.
Applying (\ref{move}) we can find an ideal $J'\subset R$ and a surjection
$\eta:R^{d-1}\sur J\cap J'$ such that $J+J'=R$, $\hh(J')\geq d-1$, and
$\eta\ot R/J=\omega_1$. We have $\hh(J')=d-1$. Let $\omega_{J'}=\eta\ot R/J'$. Then,
$(J,\omega_1)+(J',\omega_{J'})=0$ in $E^{d-1}(R)$. By the Chinese Remainder Theorem, 
$\omega_2$ and $\omega_{J'}$
will induce  $\omega_{J\cap J'}:(R/J\cap J')^{d-1}\sur J\cap J'/(J \cap J')^2$.
Then, $(J\cap J',\omega_{J\cap J'})=(J,\omega_2)+(J',\omega_{J'})$. As $\mu(J\cap J')=d-1$, by
the above theorem,   $(J\cap J',\omega_{J\cap J'})=0$, and we are done.
\qed

\smallskip
 
We answer Question \ref{ci}  in the following form.

\bt\label{murthy2}
Let $p\neq 2$ and $R$ be a smooth affine algebra of dimension $d\geq 3$ over $\k$ and $P$ be a projective $R$-module of rank $d-1$ with trivial determinant.  
Then $P\iso Q\op R$ for some $R$-module $Q$ if and only if there
is a surjection $\alpha:P\sur J$ such that $J\subset R$ is an ideal of height $d-1$ and $\mu(J)=d-1$.
\et

\proof
Follows from  (\ref{split}), (\ref{cilift}), and (\ref{zero}). 
\qed

\smallskip

\rmk
If $d\geq 5$, we can actually drop the smoothness assumption  in the above theorem. But this will come at the expense of a restrictive hypothesis on $p$, as we will soon see below.

\smallskip

We need the following variant of (\ref{cilift}). In this version we shall use the cancellation theorem of Dhorajia-Keshari \cite[3.5]{dk} instead of \cite[Theorem 7.5]{frs}, thus avoiding the restriction of smoothness. In the proof of (\ref{cilift2}) below we use some computations inside the Euler class group  $E^{d-1}(R)$. We remind the reader that in Section \ref{ECG} we only need smoothness to prove that the Euler class of a projective module is well-defined.

\bt\label{cilift2}
Let $R$ be an affine algebra over $\k$ of  dimension $d\geq 5$ such that $p$ and $(d-1)!$ are relatively prime. 
Let $J\subset R$ be an ideal of height $d-1$ such that $\mu(J)=d-1$. Suppose, it is given that $J=(b_{1},\cdots,b_{d-1})+J^2$. Then,  there exist $c_{1},\cdots,c_{d-1}\in J$
such that $J=(c_{1},\cdots,c_{d-1})$ and $b_{i}\equiv c_{i}$ mod $J^2$ for $i=1,\cdots,d-1$.
\et

\proof
Let $J=(a_1,\cdots,a_{d-1})$.
These generators will induce  $\omega:(R/J)^{d-1}\sur J/J^2$. Note that $(J,\omega)=0$ in $E^{d-1}(R)$ by definition.

Let the data $J=(b_{1},\cdots,b_{d-1})+J^2$ define $\omega':(R/J)^{d-1}\sur J/J^2$. By 
(\ref{smb}),
$(J,\omega)=(J,\ol{u}\omega)$ for some $u\in R$ which is a unit modulo $J$. Let $v\in R$ be such that 
$uv\equiv 1$ modulo $J$. Note that $(v,a_1,\cdots,a_{d-1})$ is a unimodular row over $R$ and by \cite[3.5]{dk} it is completable. 

We can follow the arguments as in \cite[Proposition, page 956, second proof]{rs1}  to conclude that
there is a matrix ${\sigma}\in M_{d-1}(R)$ with determinant ${u}^{d-2}$ modulo $J$ such that,
if $(a_1,\cdots,a_{d-1})\sigma=(f_1,\cdots,f_{d-1})$, then $J=(f_1,\cdots,f_{d-1})$. Let the generators $f_1,\cdots,f_{d-1}$ induce $\omega'':(R/J)^{d-1}\sur J/J^2$. Clearly,
$0=(J,\omega'')=(J,\ol{u}^{d-2}\omega)$.

\smallskip

\noindent
\emph{Case 1.} Let $d$ be odd. Then, $(J,\omega')=(J,\ol{u}\omega)=(J,\ol{u}^{d-2}\omega)=0$
 by
(\ref{square}). Therefore, we are done by (\ref{zero}).

\smallskip

\noindent
\emph{Case 2.} Let $d$ be even. We are required to prove that $(J,\ol{u}\omega)=0$ and for that purpose,
we may assume that $\ol{u}\omega$ is induced by $J=(ua_1,a_2,\cdots,a_{d-1})+J^2$. We may further assume by some standard arguments that $a_{d-1}$ is not a zero-divisor. Let $A=R/(a_{d-1})$ and  `tilde' denote reduction mod 
$(a_{d-1})$. It is easy to check that there is a group homomorphism $E^{d-2}(A)\lra E^{d-1}(R)$ (for an idea of proof, see \cite[7.4]{brs3}).
 Then, in $A$ we have $\widetilde{J}=(\widetilde{a_1},\cdots,\widetilde{a_{d-2}})$. 
Also,  $\widetilde{J}=(\widetilde{ua_1},\widetilde{a_2},\cdots,\widetilde{a_{d-2}})+\widetilde{J^2}$, and it is enough to lift these to generators of $\widetilde{J}$. We now apply Case 1 and the morphism
$E^{d-2}(A)\lra E^{d-1}(R)$ to conclude.
\qed

\smallskip

It is now easy to prove
the following result.

\bt\label{murthy2a}
Let $R$ be an affine algebra over $\k$ of dimension $d\geq 5$ such that $p$ and $(d-1)!$ are relatively prime. Let $P$ be a projective $R$-module of rank $d-1$ with trivial determinant.  
Then $P\iso Q\op R$ for some $R$-module $Q$ if and only if there
is a surjection $\alpha:P\sur J$ such that $J\subset R$ is an ideal of height $d-1$ and 
$\mu(J)=d-1$.
\et

\proof
By a theorem of Eisenbud-Evans \cite{ee}, there is a surjection $\beta:P\sur J$ such that $J\subset R$ is an ideal of height $d-1$.
 
Let $P\iso Q\op R$ for some $R$-module $Q$. We have proved in (\ref{mohan}) that 
$\mu(J)=d-1$.

Conversely, assume that there is a surjection $\alpha:P\sur I$  such that $I\subset R$ is an ideal of height $d-1$ and $\mu(I)=d-1$. 
Fix $\chi:R\iso \wedge^{d-1} P$. The $R/I$-module $P/IP$ is free. Choose $\sigma:(R/I)^{d-1}\iso P/IP$ such that 
$\wedge^{d-1}\sigma=\chi\ot R/I$. Let $\omega_I$ be the composite:
$$(R/I)^{d-1}\stackrel{\sigma}{\iso} P/IP\stackrel{\ol{\ga}}{\sur} I/I^2.$$
Let $\omega_I$ correspond to $I=(f_1,\cdots,f_{d-1})+I^2$. 

By the hypothesis, $\mu(I)=d-1$. We can apply (\ref{cilift2}) and obtain 
$a_1,\cdots,a_{d-1}\in I$  such that $I=(a_1,\cdots,a_{d-1})$ where $a_i-f_i\in I^2$ for $1\leq i \leq d-1$. Now we can apply (\ref{sub3}).
\qed

\section{Applications II: Projective modules over threefolds}
In this section we apply the Fasel-Rao-Swan cancellation theorem \cite[7.5]{frs} to the  results obtained in the previous sections to answer  Question \ref{quest} raised in the introduction. 
Throughout this section we assume that $p\neq 2$ and $R$ is smooth, although as remarked earlier in the introduction, the smoothness assumption may not be necessary.

\bt\label{iso}
Let $p\neq 2$, and let $R$ be a smooth affine algebra of dimension $3$ over $\k$ and $P, Q$ be  projective $R$-modules of rank 
$2$, each with trivial determinant. Fix  isomorphisms $\chi_P:R\iso \wedge^{2} P$ and 
$\chi_Q :R\iso \wedge^{2} Q$.  Then, $e_2(P,\chi_P)=e_2(Q,\chi_Q)$ in $E^2(R)$ if and only if $P$ is
isomorphic to $Q$.
\et

\proof
Let us first assume that $P$ is isomorphic to $Q$. Let $\gamma:P\iso Q$ be an isomorphism. Let $\ga:Q\sur J$ be a surjection, where $J$ is an ideal of height $2$. Suppose that $(\ga,\chi_Q)$ induces
$e_2(Q,\chi_Q)=(J,\omega)$. Note that $\wedge^2(\gamma):\wedge^2 P\iso \wedge^2 Q$ is basically 
multiplication by a unit $\ol{u}\in (R/J)^{\ast}$. It then easily follows that $(\ga\gamma,\chi_P)$ will induce $e(P,\chi_P)=(J,\ol{u}\omega)$. But $(J,\omega)=(J,\ol{u}\omega)$ by 
(\ref{equal}).

Conversely, assume that $e_2(P,\chi_P)=e_2(Q,\chi_Q)$ in $E^2(R)$. Let 
$e_2(P,\chi_P)=(J,\omega_J)=e_2(Q,\chi_Q)$. Then, by (\ref{ind}), there exist surjections $\ga:P\sur J$, 
$\gb:Q\sur J$ such that $(\ga,\chi_P)$ and $(\gb,\chi_Q)$ both induce $\omega_J$. In other words, we have a commutative diagram:

$$
\xymatrix{
        & (R/J)^2 \ar [r]^{\delta_1} \ar [d]^{\sigma} & P/JP \ar
       [r]^{\ol{\ga}}   & J/J^2  \ar [d]^{} &  \\
        & (R/J)^2 \ar [r]^{\delta_2} & Q/JQ \ar [r]^{\ol{\gb}} &
        J/J^2  &  
        }$$ 
where $\sigma\in SL_2(R/J)$, $\wedge^2\delta_1=\chi_P\ot R/J$, $\wedge^2\delta_2=\chi_Q\ot R/J$,
and the right vertical map is the identity. Let $\theta=\delta_2\sigma\delta_1^{-1}$. Then we have
$\theta:P/JP\iso Q/JQ$, $\ol{\ga}=\ol{\gb}\theta$, and 
$\wedge^2\theta=(\wedge^2\delta_2)(\wedge^2 \delta_1)^{-1}$ (as $\sigma\in SL_2(R/J)$). We can now
apply \cite[3.5]{b} to complete the proof.
\qed

\medskip

We can now answer  Question \ref{quest} raised in the introduction.

\bc
Let $p\neq 2$, and let $R$ be a smooth affine algebra of dimension $3$ over $\k$, and $P, Q$ be  projective $R$-modules of rank 
$2$, each with trivial determinant. Then, $P\iso Q$ if and only if there is an ideal $J\subset R$ of height $2$ such that both $P$ and $Q$ map onto $J$.
\ec

\proof
One way it is trivial. So let there exist an ideal $J$ of height $2$ such that there are surjections $\ga:P\sur J$ and $\gb:Q\sur J$. Fix $\chi_P:R\iso \wedge^2P$ and $\chi_Q:R\iso \wedge^2 Q$. Let $(\ga,\chi_P)$ induce $e_2(P,\chi_P)=(J,\omega_1)$ and $(\gb,\chi_Q)$ induce
$e_2(Q,\chi_Q)=(J,\omega_2)$. The proof is complete by applying (\ref{equal}) and (\ref{iso}).
\qed

\bc
Let $p\neq 2$, and let $R$ be a smooth affine algebra of dimension $3$ over $\k$ and $P$ be a projective $R$-module of rank 
$2$ with trivial determinant. Then, $P$ is cancellative.
\ec

\proof
Let $Q$ be an $R$-module of rank two such that $P\op R^n\iso Q\op R^n$ for some $n$.
As $P\op R$ is cancellative by  \cite{su}, we have $P\op R\iso Q\op R$. By \cite[6.7]{brs3},
there is an ideal $J\subset R$ of height at least two such that both $P$ and $Q$ map onto $J$.
Applying the above corollary, we have $P\iso Q$.
\qed


\begin{thebibliography}{[XXX]}
\bibitem{af1} Aravind Asok and Jean Fasel, A cohomological classification of vector bundles on smooth affine threefolds, {\it Duke Math. Journal} {\bf 163} (2014), 2561-2601.
\bibitem{af2} Aravind Asok and Jean Fasel, Splitting vector bundles outside the stable range and 
${\mathbb A}^1$-homotopy sheaves of  punctured affine spaces, {\it J. American Math. Soc.} {\bf 28} 
(2015), 1031-1062.
\bibitem{b} S. M. Bhatwadekar, Cancellation theorems for projective modules
over a two dimensional ring and its polynomial extensions, 
{\it Compositio Math.} {\bf 128} (2001), 339-359.
\bibitem{brs1} S. M. Bhatwadekar and  Raja Sridharan,
Projective generation of curves in polynomial extensions of an affine domain 
and a question of Nori, 
{\it Invent. Math.} {\bf 133} (1998), 161-192.
\bibitem{brs2} S. M. Bhatwadekar and  Raja Sridharan,
Zero cycles and the Euler class groups of smooth real affine varieties, 
{\it Invent. Math.} {\bf 136} (1999), 287-322.
\bibitem{brs3} S. M. Bhatwadekar and  Raja Sridharan, 
The Euler class group of a Noetherian ring, {\it Compositio Math.} {\bf 122}
(2000), 183-222.
\bibitem{brs4} S. M. Bhatwadekar and  Raja Sridharan,
On Euler classes and stably free projective modules,
in: Algebra, arithmetic and geometry, Part I, II(Mumbai, 2000), 139--158, 
Tata Inst. Fund. Res. Stud.
Math., 16, Tata Inst. Fund. Res., Bombay, 2002. 
\remove{
\bibitem[CHK]{chk} Jean-Louis Colliot-Th\'{e}l\`{e}ne, Raymond T. Hoobler and Bruno  Kahn, 
The Bloch-Ogus-Gabber theorem, In : Algebraic $K$-theory (Toronto, ON, 1996),  31–94,
 Fields Inst. Commun., 16, Amer. Math. Soc., Providence, RI, 1997.} 
\bibitem{c} P. M. Cohn, On the structure of the $GL_2$ of a ring, 
{\it  I.H.E.S.} {\bf 30} (1966), 5-53.
\bibitem{d} M. K. Das, The Euler class group 
of a polynomial algebra,  {\it J. Algebra} {\bf 264} (2003), 582-612.
\remove{
\bibitem[Du]{du} S. P. Dutta, A theorem on smoothness - Bass-Quillen, Chow groups and intersection multiplicity of Serre,
{\it Trans. Amer. Math. Soc.} {\bf 352} (1999), 1635-1645.}
\bibitem{dk} Alpesh M. Dhorajia and Manoj K. Keshari, note on cancellation of projective modules,
{\it J. Pure and Appl. Algebra} {\bf 216} (2012), 126-129.
\bibitem{dug} Amartya K. Dutta and Neena Gupta, The epimorphism theorem and its generalizations,
{\it Journal of Algebra and Its Applications} {\bf 14} (2015), 1540010, 30 pages.
\bibitem{ee} D. Eisenbud and E. G. Evans,  
Generating modules efficiently: Theorems from algebraic $K$-Theory, 
{\it J. Algebra} {\bf 27} (1973), 278-305.
\bibitem{frs} Jean Fasel, Ravi A. Rao and R. G. Swan, On stably free modules over affine algebras,
{\it Publ. Math. Inst. Hautes Études Sci.} {\bf 116} (2012), 223–243.
\remove{\bibitem[Ga1]{g} O. Gabber, "Some theorems on Azumaya algebras" in {\it The Brauer Group (Sem.,
Les Plans-sur-Bex, 1980)},   Lecture Notes in Math., 844, Springer, Berlin-New York, 1981, 129--209.
\bibitem[Ga2]{g2} O. Gabber, Gersten's conjecture for some complexes of vanishing cycles,
{\it Manuscripta Math.} {\bf 85} (1994), 323-343.}

\bibitem{hs} C. Huneke and I. Swanson, Integral closure of ideals, rings, and modules, 
in: London Mathematical Society Lecture Note Series, Vol. 336, Cambridge University Press, Cambridge, 2006.
\bibitem{ka} D. Katz, Generating ideals up to projective equivalence,
{\it Proc. Amer. Math. Soc.} {\bf 120} (1994), 401-414.
\bibitem{ke} Manoj K. Keshari, Euler Class group of a Noetherian ring, M.Phil. thesis, 
available at : https://arxiv.org/abs/1408.2645

\bibitem{ks} A. Krishna, V. Srinivas, Zero cycles on singular varieties, in “Algebraic Cycles and Motives I”, London
Mathematical Society Lecture note ser., Cambridge University Press, 343, (2007), 264–277.
\bibitem{la} T. Y. Lam, Serre's problem on projective modules. in: Springer Monographs in Mathematics Springer-Verlag, Berlin, 2006.

\remove{
\bibitem[Ly]{ly} G. Lyubeznik, The number of defining equations of affine algebraic sets,
{\it Amer. J. Math.} {\bf 114} (1992), 413-463.}
\bibitem{m1} S. Mandal, On efficient generation of ideals, 
{\it Invent. Math.} {\bf 75} (1984), 59-67.
\bibitem{mrs} S. Mandal and  Raja Sridharan, 
Euler classes and complete intersections, {\it J. Math. Kyoto Univ.} {\bf 36} (1996), 453-470.
\bibitem{mamu} S. Mandal and M. P. Murthy, Ideals as sections of projective modules,
{\it J. Ramanujan Math. Soc.} {\bf 13}  (1998), 51--62.
\bibitem{mo1} N. Mohan Kumar, Complete intersections, {\it J. Math. Kyoto Univ.} {\bf 17}
(1977), 533-538. 
\bibitem{mo} N. Mohan Kumar, On two conjectures about polynomial rings,
{\it Invent. Math.} {\bf 46} (1978), 225-236.
\bibitem{mo3} N. Mohan Kumar, Some theorems on generation of
ideals in affine algebras, {\it Comment.Math.Helv.} {\bf 59} (1984), 243-252.
\bibitem{mmr} N. Mohan Kumar, M. P. Murthy and Amit Roy, A cancellation theorem for projective modules over finitely generated rings, in :  Algebraic geometry and commutative algebra, Vol. I,  281--287, Kinokuniya, Tokyo, 1988.
\bibitem{mu1} M. P. Murthy, Complete intersections, in \emph{Conference on Commutative Algebra
1975} (Queen's Univ., Kingston, Ont., 1975), Queen's Papers on
Pure and Applied Math. {\bf 42}, Queen's Univ., Kingston, Ont., 1975, pp. 196-211. 
\bibitem{mu2} M. P. Murthy, A survey of obstruction theory for 
projective modules of top rank, {\it Contemp. Math.} {\bf 243}
American Mathematical Society, Providence, Rhode Island (1999).
\bibitem{mup} M. P. Murthy and C. Pedrini, $K_0$ and $K_1$ of polynomial rings, 
in: Algebraic K-Theory II, Lecture Notes in Mathematics, Springer-Verlag 342, 109-121 (1973).
\remove{\bibitem[Na]{n} B. S. Nashier, Efficient generation of ideals in polynomial rings, {\it Journal of Algebra}
{\bf 85} (1983), 287-302.
\bibitem{ni} Y. Nisnevich, Stratified canonical forms of matrix valued functions in a neighborhood of a transition point,
{\it Internat. Math. Res. Notices.} {\bf 10} (1998), 513-527.  
\bibitem{p} D.  Popescu, Polynomial rings and their projective modules,
{\it Nagoya Math. J.} {\bf 113} (1989), 121-128.}
\bibitem{q} D. Quillen, Projective modules over polynomial rings, 
{\it Invent. Math.} {\bf 36} (1976), 167-171.
\bibitem{ra} R. A. Rao, An elementary transformation of a special unimodular vector
to its top coefficient vector, {\it Proc. Amer. Math Soc.} {\bf 93} (1985), 21-24.
\bibitem{ra2} R. A. Rao, The Bass-Quillen conjecture in dimension three but
characteristic $\neq 2$,$3$ via a question of A. Suslin, {\it Invent. Math.}
{\bf 93} (1988), 609-618.
\bibitem {ra3} R. A. Rao, On completing unimodular polynomial vectors of length three,
{\it Trans. Amer. Math. Soc.} {\bf 325} (1991), 231-239.
\bibitem{se} J.-P. Serre, Modules projectifs et espaces fibres a fibre 
vectorielle, {\it Sem. Dubreil-Pisot} {\bf 23} (1957/58).
\bibitem{rs1} Raja Sridharan, Non-vanishing sections of algebraic
vector bundles,{\it J. Algebra} {\bf 176} (1995), 947-958.
\bibitem{rs} Raja Sridharan, Projective modules and complete 
intersections, {\it K-Theory} {\bf 13} (1998), 269-278.
\bibitem{sr} V. Srinivas, $K_1$ of the cone over a curve,
J. reine angew. Math. {\bf 381} (1987), 37—50.
\bibitem{su} A. A. Suslin, A cancellation theorem for projective modules over algebras, 
{\it Soviet.
Math. Dokl.} {\bf 18} (1977), 1281-1284.
\bibitem{sv} A. A. Suslin and L. N. Vaserstein, Serre's problem on projective modules over polynomial
rings and algebraic $K$-theory, {\it Izv. Akad. Nauk. SSSR, Ser. Mat.} {\bf 40} (1976), 937-1001.
\remove{\bibitem[Sw]{sw} R. G. Swan,  N\'{e}ron-Popescu desingularization, in : Algebra and Geometry (Taipei, 1995),
Lect. Algebra Geom. vol 2, International Press, Cambridge, MA, 135-192.}
\bibitem{sw} Richard G. Swan, Serre's problem. \emph{Conference on Commutative Algebra 1975} 
(Queen's Univ., Kingston, Ont., 1975), pp. 1–60. Queen's Papers on Pure and Applied Math., {\bf 42}, Queen's Univ., Kingston, Ont., 1975.
\bibitem{v} L. N. Vaserstein, On the stabilization of the general linear group over
a ring. {\it Mat. Sb. (N.S)} 79 (121) 405–424 (Russian); translated as Math.
USSR-Sb.8, 1969, 383–400.
\end{thebibliography}
\end{document}